\documentclass[12pt]{amsart}
\usepackage{amscd, amssymb, amsfonts, amsthm, mathrsfs}
\usepackage{fullpage}
\usepackage{quiver}
\usepackage{dynkin-diagrams}
\usepackage[hidelinks]{hyperref}
\usepackage[shortlabels]{enumitem}

\newtheorem{theorem}{Theorem}[section]
\newtheorem{lemma}[theorem]{Lemma}
\newtheorem{proposition}[theorem]{Proposition}
\newtheorem{corollary}[theorem]{Corollary}

\theoremstyle{definition}
\newtheorem{remark}[theorem]{Remark}
\newtheorem{example}[theorem]{Example}
\newtheorem{definition}[theorem]{Definition}

\title{Buildings governed by twisted Bruhat orders}

\author[H.\ Sykes]{Harvey Sykes}
\address{Fakult\"at f\"ur Mathematik,
Ruhr-Universit\"at Bochum,
D-44801 Bochum, Germany}
\email{harvey.sykes@rub.de}

\begin{document}

\begin{abstract}
    This paper introduces the notion of a twisted building associated to a Coxeter system and an initial section of a reflection order. These are pairs of buildings and a map between them which depends on a twisted Bruhat order. It is proved that projection maps exist between components and that the common apartment property holds for a suitable notion of an apartment, generalising what is known for twin buildings. Several well-known examples are given, including a family of twisted buildings associated to an affine building and its fa\c{c}ades.
\end{abstract}

\maketitle

\tableofcontents

\section{Introduction}
The notion of a building was introduced by Tits in the 1950s with the aim to construct a space upon which a complex simple Lie group acts; specifically, those of exceptional type. In this context, the space is built from the parabolic subgroups and behaves locally with respect to the finite Weyl group. Examples of this kind are called spherical buildings and have numerous applications, such as Serre's notion of $G$-complete reducibility \cite{serre05}. In general, one can associate a building to any Coxeter system. Most, but not all, are realised using parabolic-like subgroups of an appropriate group, e.g.~parahoric subgroups for Bruhat-Tits buildings \cite{bruhat-tits84}. Since their inception, the definition has been recast several times, each with a different emphasis, be it apartments \cite{tits74}, chambers \cite{tits81}, or metrics \cite{davis98}. A more recent formulation is the $W$-metric approach \cite{tits92}, which focusses on a Weyl-distance map 
\begin{equation*}
    \delta:X\times X\rightarrow W
\end{equation*}
satisfying some natural axioms with respect to a standard length function on the Coxeter group $W$. These axioms allow one to define projection maps quickly, while apartments become trickier to establish, cf.~\cite[\S5]{abramenko-brown08}.

In the 1980s, twin buildings were defined and studied by Ronan and Tits. As explained in \cite{tits92}, one motivation for this work was the emerging theory of Kac-Moody groups. Roughly speaking, a Kac-Moody group is to an arbitrary Weyl group what a simple algebraic group is to a finite Weyl group. For infinite Weyl groups, the corresponding Kac-Moody group admits two classes of parabolic subgroup. This reflects the fact that the ambient root system admits two distinct conjugacy classes of standard positive systems. A twin building encodes these conjugacy classes with two buildings $X_+$ and $X_-$ of the same type, along with a map
\begin{equation*}
    \delta^*:(X_+\times X_-)\cup(X_-\times X_+)\rightarrow W
\end{equation*}
satisfying a similar set of axioms as before, after replacing the standard length function with its negation. In particular, the buildings are governed by a length function opposite to the one controlling how they relate.

This article introduces a new object called a \emph{twisted building} (see \S\ref{twisted:def}) which generalises the two examples above. A twisted building comprises a building $X$ and a satellite building $X_A$ whose type is a Coxeter system $(W_A,S_A)$ associated to an initial section $A$, cf.~\S\ref{cox:initial-sections}. The buildings are related by a map
\begin{equation*}
    \delta^*:(X\times X_A)\cup(X_A\times X)\rightarrow W
\end{equation*}
which satisfies the same axioms again, this time replacing the standard (resp.~opposite) length function with Dyer's twisted length function \cite{dyer-heckeii92}. Both the length function and the type of $X_A$ depend on $A$. At the extremes, when $A$ is empty or is the set of all reflections, one retrieves the definition of a building and a twin building, respectively. New examples appear when $X$ is not spherical. In particular, if $X$ is an affine building and $X_A$ is a fa\c{c}ade of $X$, then a twisted building can be defined using the chimney retractions of \cite{milicevic-gallery24}. This construction is made explicit in \S\ref{aff}.

Twisted buildings admit projection maps just like ordinary and twin buildings. This is proved by adapting well-known arguments, first for panels (Lemma \ref{twisted:proj:lem:panels}) and then to all spherical residues (Proposition \ref{twisted:proj:prop:residues}). In buildings, projections exist no matter the type of the target residue. By contrast, projections between components of twin buildings are well-defined only if the target is spherical. The twisted buildings lie somewhere in the middle. They may, but not necessarily, admit projections onto non-spherical residues. A precise characterisation for existence is given with respect to the ambient initial section (Proposition \ref{twisted:proj:prop:existence}). Remarkably, this sometimes allows for projections onto the whole of a non-spherical satellite building (see Corollary \ref{aff:twist:cor:proj}). The analogue of a twin apartment is defined in \S\ref{twisted:apartments}. It is proved that apartments are convex in a suitable sense (Lemma \ref{twisted:apartment:lem:convex}) and can be produced by extending a building isometry in $X$ (Theorem \ref{twisted:apartments:thm}). This implies that any two chambers that lie in different components share a common twisted apartment (Corollary \ref{twisted:apartments:cor:common}). The same is not true for arbitrary pairs of chambers (see Example \ref{twisted:apartments:ex}).

The affine examples in \S\ref{aff} show how twisted buildings can be used to describe the relationship between an affine building and its fa\c{c}ades. Each admits a building structure, and their chambers interact through a twisting map (Theorem \ref{aff:twist:thm}). This gives a conceptual reason for the behaviour of chimney germs (see Proposition \ref{aff:fac:prop:chimneys}) and the existence of vectorial projection maps from the affine building to its fa\c{c}ades (Corollary \ref{aff:twist:cor:proj}). In particular, projections exist onto an entire non-spherical fa\c{c}ade. There is also a group-theoretic interpretation to these constructions.  When $X$ is the Bruhat-Tits building of a reductive group over a local field, the fa\c{c}ade $X_F$ can be identified with the (reduced) Bruhat-Tits building of a Levi factor. In the split case, the twisting map can be deduced from certain Bruhat-type decompositions, cf.~\cite[\S11.2]{goertz10} or \cite[\S5]{milicevic-gallery24}. This will not be developed further in this article.

The motivation for this work is from Kac-Moody groups. In Weyl groups, there is an initial section of the form $A=A_\Psi$ for every (possibly non-standard) positive system $\Psi$ of the root system. The standard positive and negative systems give the two extremes. Less prevalent in the literature are the non-standard systems, which provide interesting twisted buildings that are neither buildings nor twin buildings (see \S\ref{aff}). Furthermore, there are subgroups in Kac-Moody groups that realise them: the \emph{quasi-parabolic} subgroups of \cite{sykes-thesis26}. These are not parabolic in the usual sense but exhibit similar properties, such as a Bruhat-type decomposition \cite{billig-dyer94}. More details will appear in forthcoming work.

The layout of this article is as follows. In \S\ref{cox}, the theory surrounding initial sections is recalled. This includes twisted length functions and reflection subgroups. The section ends with unpublished results of Dyer on minimal twisted length representatives in cosets of parabolic subgroups. The central definition is given in \S\ref{twisted}. This is where the basic properties of twisted buildings, such as projection maps and apartments, are proved. The article finishes with \S\ref{aff}. In this section, a family of twisted buildings is constructed from an affine building and its fa\c{c}ades.

\begingroup
\makeatletter
\let\@tocwrite\@gobbletwo
\subsection*{Acknowledgments}
\leavevmode
\endgroup
The majority of this work appears, at least implicitly, within the author's PhD thesis from the University of Essex. Throughout this period, the author was funded by a Department of Mathematical Sciences Doctoral Scholarship. The author would like to thank Alastair Litterick for their supervision and guidance during this time. The author is also grateful to Bernhard M\"{u}hlherr for their interest and feedback on an earlier manuscript. Finally, the author wishes to thank Matthew Dyer for their generosity in sharing the proofs that appear in Appendix \ref{ap:proofs}.

\section{Preliminaries}\label{cox}
This section records the prerequisites for the rest of the paper and contains some extensions of standard facts in Coxeter theory. The new results presented below are either implicit in the literature or currently unpublished. The proofs were communicated to the author privately by Dyer, and are reproduced in Appendix \ref{ap:proofs} for completeness.

\subsection{Notation}
Throughout, let $(W,S)$ be a Coxeter system with root system $\Delta$. Denote by $\Delta^+$ the standard positive system of $\Delta$. In particular, there is a canonical bijection $\alpha\mapsto r_\alpha$ from the standard positive system $\Delta^+$ to the set of reflections $T$ of $W$. Set $\mathbb{N}=\{1,2,\dotsb\}$.

\subsection{Twisted length functions}\label{cox:initial-sections}
The standard length function is an essential feature of a Coxeter system. This has been generalised by Dyer and plays a role in Kazhdan-Lusztig theory. Unless stated otherwise, the following material can be found in \cite{dyer-hecke93,dyer-heckeii92}. 

A \emph{reflection order} $\preccurlyeq$ on $T$ is a total order satisfying
\begin{equation*}
    r_\alpha\preccurlyeq r_\gamma\implies r_\alpha\preccurlyeq r_\beta\preccurlyeq r_\gamma
\end{equation*}
for all $\alpha,\beta,\gamma\in\Delta^+$ with $\beta=c\alpha+d\gamma$ for some $c,d\geq0$.
An \emph{initial section} of the reflection order $\preccurlyeq$ is a subset $A\subseteq T$ such that $r\prec s$ for all $r\in A$ and $s\in A^c$. Denote by $\mathcal{A}$ the initial sections of reflection orders on $T$. There may be several reflection orders associated to an initial section, so elements of $\mathcal{A}$ will just be called initial sections.

\begin{example}\label{cox:initial-sections:ex}
    A natural class of initial sections arises from \emph{positive systems} of roots, i.e.~closed subsets $\Psi$ of $\Delta$ satisfying
    \begin{equation*}
        \Psi\cap(-\Psi)=\emptyset\quad\text{and}\quad\Psi\cup(-\Psi)=\Delta.
    \end{equation*}
    These were defined in \cite{jakobsen-kac83} where it was proved that there are finitely many distinct $W$-orbits of positive systems when $W$ is affine, enumerated (up to a sign) by the subsets of simple roots of the underlying finite root system. In general, the number of orbits of positive systems may be infinite. For every positive system $\Psi$ of $\Delta$, the set 
    \begin{equation*}
    A_\Psi=\{r_\alpha\in T\ |\ \alpha\in\Delta^+\cap(-\Psi)\}
\end{equation*}
    is an initial section of a reflection order on $T$, cf.~\cite[\S1.2]{dyer-quot94}.
\end{example}

Write $\mathcal{P}(X)$ for the power set of $X$. There is a function $N:W\rightarrow \mathcal{P}(T)$ characterised by $N(s)=\{s\}$ for all $s\in S$ and $N(vw)=N(v)+vN(w)v^{-1}$ for all $v,w\in W$. Here and elsewhere, the addition sign refers to the operation of symmetric difference. This induces an action of $W$ on the subsets $A\subseteq T$, where $wA=N(w)+wAw^{-1}$ for all $w\in W$. The next result is \cite[Lem.~2.7]{dyer-hecke93}.

\begin{lemma}\label{cox:twisted-length:lem:init}
    Let $A$ be an initial section and $w\in W$. Then $A^c$ and $wA$ are initial sections.
\end{lemma}

Let $A$ be an initial section. Define the \emph{twisted length function} $\ell_A:W\rightarrow\mathbb{Z}$ by
\begin{align*}
    \ell_A(w)&=\ell(w)-2|N(w^{-1})\cap A|.
\end{align*}
That is to say, $\ell_A(w)$ counts the number of roots in $\Delta^+$ sent to $-\Delta^+$ by $w$, where a root is counted with multiplicity $-1$ if its corresponding reflection lies in $A$ and with multiplicity $1$ otherwise. The standard length function $\ell$ and its opposite $-\ell$ are obtained when $A$ is empty or the set of all reflections, respectively. The following proposition records the basic properties of twisted length functions.

\begin{proposition}\label{cox:length:prop}
    Let $A$ be an initial section. Then, for all $s\in S$, $t\in T$, and $w\in W$, the following are true:
    \begin{enumerate}[series=twisted-lengths]
        \item\label{cox:length:prop:identity} $\ell_A(e)=0$.
        \item\label{cox:length:prop:negate} $\ell_{A^c}(w)=-\ell_A(w)$.
        \item\label{cox:length:prop:mult-char}
        \begin{enumerate}
            \item $\ell_A(tw)>\ell_A(w)$ if and only if $t\notin wA$.
            \item $\ell_A(wt)>\ell_A(w)$ if and only if $t\notin N(w^{-1})+A$.
        \end{enumerate}
        \item\label{cox:length:prop:mult-simple} $|\ell_A(sw)-\ell_A(w)|=1$.
    \end{enumerate}
\end{proposition}

Every initial section $A$ defines a \emph{twisted Bruhat order} $\leq_A$ on $W$ such that 
\begin{equation*}
    v\leq_A w\iff \text{there is}\ v=v_0,v_1,\dots,v_n=w\ \text{with}\ \ell_A(v_i)=\ell_A(v_{i-1})+1\ \text{and}\ v_iv_{i-1}^{-1}\in T.
\end{equation*}
Similar to before, the Bruhat order and its reverse are obtained when $A$ is empty or equal to $T$. There is also an analogue of Lemma \ref{cox:twisted-length:lem:init} which will be used without reference, cf.~\cite[1.12]{dyer-quot94}.

\subsection{Reflection subgroups}\label{cox:refl-subgps}
To every initial section one can associate a Coxeter system which interacts well with the corresponding twisted length function. The main reference for this construction is \cite{dyer-heckeii92}.

A \emph{reflection subgroup} of $W$ is a subgroup generated by reflections. Any initial section $A$ defines a reflection subgroup $W_A$ of $W$ generated by 
\begin{equation*}
    S_A=\{t\in T\ |\ tA=\{t\}+A\}.
\end{equation*}
The pair $(W_A,S_A)$ is a Coxeter system, but $S_A$ need not be the canonical generating set of Deodhar and Dyer \cite{dyer-refl90,deodhar89}. The corresponding map $N_A:W_A\rightarrow\mathcal{P}(T\cap W_A)$ satisfies $N_A(w)=A+wA$ for all $w\in W_A$. Denote by $T_A=T\cap W_A$ the reflections of $W_A$. When $A$ is empty or equal to $T$, one has $S=S_A$ and $W=W_A$. In general, $W_A$ may be proper and in some cases trivial. The next result can be found in \cite[\S~1.8]{dyer-heckeii92}.

\begin{proposition}\label{cox:refl-subgps:prop}
    Let $A$ be an initial section. Then 
    \begin{equation*}
        |\ell_A(wt)-\ell_A(w)|=1
    \end{equation*}
    for all $t\in S_A$ and $w\in W$.
\end{proposition}

Before continuing some notation must be set. Denote by $\ell'$ the standard length function on $(W_A,S_A)$. This need not coincide with the restriction of $\ell$ to $W_A$. The intersection $A\cap T_A$ is also an initial section for the Coxeter system $(W_A,S_A)$. As such, one defines a third length function by
\begin{equation*}
    \ell_A'(w)=\ell'(w)-2|N_A(w^{-1})\cap A|
\end{equation*}
for all $w\in W_A$.

\subsection{Minimal length representatives}\label{cox:cosets}
The theory of projection maps in buildings is closely related to elements of minimal length in parabolic subgroups and their cosets. Analogously, projection maps defined in twisted buildings will be linked to parabolic subgroups of reflection subgroups, and elements of minimal twisted length. The unpublished results that follow are due to Dyer in their full generality, cf.~Appendix \ref{ap:proofs}.

Let $A$ be an initial section. Set 
\begin{align*}
    X_{W_A}&=\{w\in W\ |\ \ell(tw)>\ell(w)\ \text{for all}\ t\in T_A\},\\
    &=\{w\in W\ |\ N(w)\cap T_A=\emptyset\}.
\end{align*}
The next proposition extends a well-known result for standard parabolic subgroups to arbitrary reflection subgroups. For the general case, see \cite[Prop.~1.13]{dyer1987hecke} and \cite[Cor.~3.4]{dyer-refl90}.
\begin{proposition}\label{cox:cosets:prop:X_W}
    The set $X_{W_A}$ consists of the elements of minimal length in their right coset of $W_A$. Moreover, every element $w$ of $W$ can be uniquely expressed as $w=uv$ where $u\in W_A$ and $v\in X_{W_A}$.
\end{proposition}

In particular, for all $w\in W$ there is a unique element $x_0\in X_{W_A}$ contained in $W_Aw$ with minimal length. There is a left analogue of Proposition \ref{cox:cosets:prop:X_W} obtained by inverting elements which will be used freely in the sequel.

Let $J\subseteq S$ and $K\subseteq S_A$. Set
\begin{align*}
    ^JW_A&=\{w\in W\ |\ \ell_A(sw)>\ell_A(w)\ \text{for all}\ s\in J\},\\
    &=\{w\in W\ |\ W_J\cap wA=\emptyset\},\ \text{and}\\
    W_A^K&=\{w\in W\ |\ \ell_A(wt)>\ell_A(w)\ \text{for all}\ t\in K\},\\
    &=\{w\in W\ |\ (N(w^{-1})+A)\cap (W_A)_K=\emptyset\}.
\end{align*}
When $A$ is empty, write ${^JW}$ and $W^K$ respectively. The two gate properties that follow anticipate what will later be proved for projections between the components of a twisted building. For the proof, see Appendix \ref{ap:proofs}.

\begin{proposition}\label{cox:cosets:prop:posets}
    \begin{enumerate}
        \item\label{cox:cosets:prop:posets:left} Let $w=uv$ where $u\in W_J$ and $v\in{^JW_A}$. Then
        \begin{equation*}
            \ell_A(w)=\ell(u)+\ell_A(v).
        \end{equation*}
        In particular, the map $W_J\rightarrow W_Jv$ defined by $u\mapsto uv$ is an isomorphism of posets, where $W_J$ has the Bruhat order induced from $W$ and $W_Jv$ has the order $\leq_A$ induced from $W$.
        \item\label{cox:cosets:prop:posets:right} Let $w=uv$ where $u\in W_A^K$ and $v\in (W_A)_K$. Then
        \begin{equation*}
            \ell_A(w)=\ell_A(u)+\ell'(v).
        \end{equation*}
        In particular, the map $(W_A)_K\rightarrow u(W_A)_K$ defined by $v^{-1}\mapsto uv$ is an isomorphism of posets, where $(W_A)_K$ has the Bruhat order induced from $W_A$ and $u(W_A)_K$ has the order $\leq_A$ induced from $W$.
    \end{enumerate}
\end{proposition}

It follows that ${^JW_A}$ is the collection of elements in $W$ of minimal twisted length $\ell_A$ in their right coset of $W_J$. Similarly, the elements of $W_A^K$ are those in $W$ of minimal twisted length $\ell_A$ in their left coset of $(W_A)_K$. If $W_J$ is finite or $A$ is empty, then every element of $W$ can be uniquely expressed as $w=uv$ where $u\in W_J$ and $v\in{^JW_A}$. Therefore, every coset $W_Jw$ has a representative of minimal twisted length $\ell_A$ which lies in ${^JW_A}$. The analogous statement is true for the cosets $w(W_A)_K$ when $(W_A)_K$ is finite or $A$ is empty. Another useful consequence of Proposition \ref{cox:cosets:prop:posets} is the following result which is likely known to experts.

\begin{corollary}\label{cox:cosets:cor:lengths}
    Let $w\in W_A$. Then
    \begin{equation*}
        \ell_A(w)=\ell'(w)-2|N_A(w)\cap A|.
    \end{equation*}
\end{corollary}
\begin{proof}
    Applying Lemma \ref{cox:cosets:lem:poset-isom-left} to $u=1$ and $v=w\in W_A$ gives $\ell_A(w)=\ell'_A(w^{-1})$. Expanding the latter expression, one obtains
    \begin{align*}
        \ell_A(w)&=\ell'(w^{-1})-2|N_A(w)\cap A|,\\
        &=\ell'(w)-2|N_A(w)\cap A|.
    \end{align*}
\end{proof}

More generally, one can study minimal length representatives in double cosets. Set
\begin{equation*}
    ^JW_A^K={^JW_A}\cap{W_A^K}.
\end{equation*}
As before, if $A$ is empty then suppress the subscript. The next result will be used to describe projections of residues onto residues in twisted buildings.

\begin{proposition}\label{ctd:prelim:cosets:prop:kilmoyer}
    Let $J\subseteq S$ and $K\subseteq S_A$. Suppose that $w\in{^JW_A^K}$. Then
    \begin{equation*}
        W_J\cap w(W_A)_Kw^{-1}=W_{L}
    \end{equation*}
    where $L={J\cap wKw^{-1}}$.
\end{proposition}

When $A$ is empty one retrieves a well-known result of Kilmoyer, cf.~\cite[Prop.~2.23]{abramenko-brown08}.

\begin{proposition}\label{ctd:prelim:cosets:prop:howlett}
    Let $J\subseteq S$ and $K\subseteq S_A$. Suppose that $v\in{^JW_A^K}$. Then every element $w$ of $W_Jv(W_A)_K$ can be uniquely expressed as
    \begin{equation*}
        w=abvc
    \end{equation*}
    where $a\in W_J\cap W^{J\cap vKv^{-1}}$, $b\in W_{J\cap vKv^{-1}}$, and $c\in (W_A)_K\cap{^{K\cap v^{-1}Jv}}(W_A)_\emptyset$. Moreover, one has 
    \begin{equation*}
        \ell_A(w)=\ell(a)+\ell(b)+\ell_A(v)+\ell'(c)
    \end{equation*}
    with $vc\in{^JW_A}$, $av\in W_A^K$, and $w=avb'c$ where $b'=v^{-1}bv\in (W_A)_{K\cap v^{-1}Jv}$ and $\ell(b)=\ell'(b')$.
\end{proposition}

In particular, if $W_J$ and $(W_A)_K$ are finite or $A$ is empty, then every double coset $W_Jv(W_A)_K$ has a unique element $w$ of minimal twisted length $\ell_A$ and $w\in{^JW_A^K}$. The latter case is a well-known result of Howlett, cf.~\cite[Cor.~2.25]{abramenko-brown08}.

\subsection{Buildings}\label{cox:buildings}
As mentioned in the introduction, this article follows the $W$-metric approach to buildings. This will mean that projection maps can be defined straight away, and characterised in detail using the results of \S\ref{cox:cosets}. More detail can be found in \cite[\S5]{abramenko-brown08}.

A \emph{building} of type $(W,S)$ is a pair $(X,\delta)$ consisting of a set $X$ of \emph{chambers}, equipped with a $W$-valued map 
\begin{equation*}
    \delta:X\times X\rightarrow W
\end{equation*}
such that, for all $x,y\in X$ with $\delta(x,y)=w$, and all $s\in S$, the following properties hold:
\begin{enumerate}[label=\bf(B\arabic*), ref=\bf(B\arabic*)]
    \item $\delta(x,y)=1$ if and only if $x=y$.
    \item If $z\in X$ satisfies $\delta(z,x)=s$ then $\delta(z,y)\in\{w,sw\}$. Furthermore, if $\ell(sw)>\ell(w)$ then $\delta(z,y)=sw$.
    \item There exists $z\in X$ with $\delta(z,x)=s$ and $\delta(z,y)=sw$.
\end{enumerate}
If $(W,S)$ is finite (resp.~affine) then $(X,\delta)$ is called \emph{spherical} (resp.~\emph{affine}). Let $(X,\delta)$ be a building of type $(W,S)$. A pair of chambers $x$ and $y$ in $X$ are said to be \emph{$s$-adjacent} if $\delta(x,y)=s$ and $s\in S$. This will be abbreviated to just \emph{adjacent} when $s$ is not specified. A \emph{gallery} from $x$ to $y$ is a sequence of pairwise consecutive adjacent chambers. The \emph{gallery distance} between two chambers is the length of a minimal gallery connecting them, which always exists. Equivalently, it is defined to be
\begin{equation*}
    d(x,y)=\ell(\delta(x,y))
\end{equation*}
for all $x,y\in X$. 

Let $J\subseteq S$. A \emph{residue} $R$ in $X$ of type $J$, or \emph{$J$-residue}, is a set of the form
\begin{equation*}
    R=\{y\in X\ |\ \delta(x,y)\in W_J\}
\end{equation*}
for some $x\in X$. Importantly, a residue of type $J$ is a building of type $(W_J,J)$. A \emph{panel} $P$ in $X$ is a residue of type $s$ for some $s\in S$. More specifically, $P$ is an \emph{$s$-panel}. 

\begin{proposition}\label{cox:buildings:prop:proj}
    Let $R$ be a residue of $X$. Then, for any chamber $x\in X$, there is a unique chamber $p$ in $R$ whose distance to $x$ is minimal. Moreover, one has
    \begin{equation*}
        \delta(x,y)=\delta(x,p)\delta(p,y)\quad\text{and}\quad d(x,y)=d(x,p)+d(p,y).
    \end{equation*} 
    for all $y\in R$.
\end{proposition}

The unique chamber of Proposition \ref{cox:buildings:prop:proj} is called the \emph{projection} of $x$ onto $R$. It will often be denoted by $\operatorname{proj}_R(x)$. The decompositions of $\delta$ and $d$ with respect to the projection are instances of the \emph{gate property}. It is also possible to define projection maps between arbitrary residues in $X$. For more on this, see \cite[\S5.3.2]{abramenko-brown08}.

A subset $Y$ of $X$ is \emph{convex} if it is closed with respect to minimal galleries. For example, a residue is convex. A subset $Y$ of $X$ is \emph{thin} if $Y\cap P$ consists of two chambers for every panel $P$ that meets $Y$. The \emph{canonical thin building} of type $(W,S)$ is $W$ itself, equipped with the map
\begin{equation*}
    \delta_W:W\times W\rightarrow W,\ \delta_W(v,w)=v^{-1}w.
\end{equation*}
An \emph{apartment} $\Sigma$ is a non-empty subset of $X$ that is convex and thin. The apartments provide a local model for $X$. Indeed, every pair of chambers has an apartment containing them. Moreover, every apartment of $X$ is isomorphic to the canonical thin building on $W$.

\section{Twisted buildings}\label{twisted}
This section defines the central object of this paper and establishes the existence of projection maps and apartments. The progression of results is inspired by the treatment of twin buildings in \cite{abramenko-brown08}.

\subsection{Definition}\label{twisted:def}
Let $A$ be an initial section.

\begin{definition}
    A \emph{twisted building} of type $(W,S,A)$ is a pair of buildings $(X,\delta)$ and $(X_A,\delta_A)$ of type $(W,S)$ and type $(W_A,S_A)$ respectively, along with a \emph{twisting} of type $A$
    \begin{equation*}
        \delta^*:(X\times X_A)\cup(X_A\times X)\rightarrow W
    \end{equation*}
    such that, for all $x\in X$ and $x_A\in X_A$ with $\delta^*(x,x_A)=w$, and all $s\in S$ and $t\in S_A$, the following properties hold:
    \begin{enumerate}[label=\bf(T\arabic*), ref=\bf(T\arabic*)]
        \item\label{twisted:def:def:1} $\delta^*(x,x_A)=\delta^*(x_A,x)^{-1}$.
        \item\label{twisted:def:def:2}
        \begin{enumerate}
            \item\label{twisted:def:def:2a} If $y\in X$ satisfies $\delta(y,x)=s$ and $\ell_A(sw)>\ell_A(w)$, then $\delta^*(y,x_A)=sw$.
            \item\label{twisted:def:def:2b} If $y_A\in X_A$ satisfies $\delta_A(x_A,y_A)=t$ and $\ell_A(wt)>\ell_A(w)$, then $\delta^*(x,y_A)=wt$.
        \end{enumerate}
        \item\label{twisted:def:def:3}
        \begin{enumerate}
            \item\label{twisted:def:def:3a} There exists $y\in X$ with $\delta(y,x)=s$ and $\delta^*(y,x_A)=sw$.
            \item\label{twisted:def:def:3b} There exists $y_A\in X_A$ with $\delta_A(x_A,y_A)=t$ and $\delta^*(x,y_A)=wt$.
        \end{enumerate}
    \end{enumerate}
\end{definition}

To the author's knowledge, twisted buildings are not defined in the literature in this level of generality. However, there are many well-known special cases of twisted buildings. 

\begin{example}\label{twisted:def:ex-standard}
    \begin{enumerate}
        \item A twisted building of type $(W,S,\emptyset)$ is a building of type $(W,S)$.
        \item A twisted building of type $(W,S,T)$ is a twin building of type $(W,S)$.
        \item\label{twisted:def:ex-standard3} Let $(X,\delta)$ be a building of type $(W,S)$ and $J\subseteq S$ spherical. Denote by $w_J$ the longest element in $W_J$. Suppose that $X_A$ is a copy of $X$, and fix a bijection sending $x_A$ to $x$. Define the Weyl-distance
        \begin{equation*}
            \delta_A(x_A,y_A)=w_J\delta(x,y)w_J
        \end{equation*}
        for all $x_A,y_A\in X_A$. Then the translate
        \begin{equation*}
            \delta^*(x,y_A)=\delta(x,y)w_J
        \end{equation*}
        is a twisting of type $T\cap W_J$ between $X$ and $X_A$. A twisting of type $T\setminus(T\cap W_J)$ can be constructed in a similar way from a twin building.
    \end{enumerate}
\end{example}

Note that the Coxeter systems associated to an initial section and its complement are the same. As a result, a pair of buildings may admit distinct twisting maps of type $A$ and of type $A^c$. This is the difference between a building and a twin building. Another similarity is the fact that a twisted building can be produced from any Coxeter group and reflection subgroup arising from an initial section.

\begin{example}\label{twisted:def:ex-thin}
    Recall that the Coxeter groups $W$ and $W_A$ are canonically thin buildings of type $(W,S)$ and $(W_A,S_A)$ respectively. Furthermore, the map $\delta^*:W\times W_A\rightarrow W$ defined by 
    \begin{equation*}
        \delta^*(v,w)=v^{-1}w
    \end{equation*} 
    for all $v\in W$ and $w\in W_A$ is a twisting. In fact, $\delta^*$ is a twisting of type $A$ and of type $A^c$. This twisted building is called the \emph{thin twisted building} associated to $W$ and $W_A$.
\end{example}

A more interesting example arises in the context of affine buildings, offering a different perspective for the notion of a retraction from infinity as defined in \cite[p.~171]{brown-book89}.

\begin{example}\label{twisted:def:ex-retract}
    Let $(X,\delta)$ be an affine building of type $(W,S)$. Consider the initial section
    \begin{equation*}
        A=\{r_{-\alpha+n\delta}\ |\ \alpha\in\Phi^+,\ n\in\mathbb{N}\}
    \end{equation*}
    where $\Phi^+$ is the standard positive system of the underlying finite root system of $W$. This is the initial section of Example \ref{cox:initial-sections:ex} associated to the so-called \emph{natural} positive system \cite[p.~7]{jakobsen-kac83}. Then $S_A$ is empty and $W_A$ is trivial. Therefore, a building of type $(W_A,S_A)$ is a point. The twisted Bruhat order $\leq_A$ is the semi-infinite order \cite[Prop.~1.16]{dyer-quot94} and the retraction from infinity
    \begin{equation*}
        r_\infty: X\rightarrow W
    \end{equation*}
    induces a twisting map of type $A$. Similarly, the retraction from minus infinity defines a twisting map of type $A^c$. This will be extended to satellite buildings of higher rank in \S\ref{aff}.
\end{example}

\begin{remark}   
    A finite rank Coxeter system may admit a reflection subsystem of infinite rank, e.g.~\cite[Ex.~3.19]{dyer1987hecke}. It is not clear to the author whether such examples arise from initial sections. In particular, there may exist twisted buildings whose components have different finitude. 
\end{remark}

For the rest of the section, fix an initial section $A$ and a twisted building $(X,X_A,\delta^*)$ of type $(W,S,C)$ where $C\in\{A,A^c\}$.

\subsection{Projections}\label{twisted:proj}
\begin{definition}
    Let $x\in X$ and $x_A\in X_A$. The \emph{distance} between $x$ and $x_A$ is
    \begin{equation*}
        d^C(x,x_A)=\ell_C(\delta^*(x,x_A)).
    \end{equation*}
\end{definition}

Note that the distance in a twisted building does not generalise the notion of codistance in a twin building. For a twin building, the distance defined above is the negation of the codistance. The convention taken in this article is justified by the following definition.

\begin{definition}
Let $P$ be a panel of $X$ and $P_A$ a panel of $X_A$. Then:
    \begin{itemize}
        \item The \emph{projection} of a chamber $x_A\in X_A$ onto $P$ is the unique chamber $\operatorname{proj}_P(x_A)$ of $P$ that minimises $d^C(-,x_A)$.
        \item The \emph{projection} of a chamber $x\in X$ onto $P_A$ is the unique chamber $\operatorname{proj}_{P_A}(x)$ of $P_A$ that minimises $d^C(x,-)$.
    \end{itemize}
\end{definition}

The next lemma proves the existence and uniqueness of panel projections.

\begin{lemma}\label{twisted:proj:lem:panels}
    Let $x\in X$ and $x_A\in X_A$. Suppose that $\delta^*(x,x_A)=w$. Then, for all $s\in S$ and $t\in S_A$, the following are true:
    \begin{enumerate}
        \item\label{twisted:proj:lem:panels:1} If $y\in X$ with $\delta(y,x)=s$ then $\delta^*(y,x_A)\in\{w,sw\}$.
        \item\label{twisted:proj:lem:panels:2} If $y_A\in X_A$ with $\delta_A(x_A,y_A)=t$ then $\delta^*(x,y_A)\in\{w,wt\}$.
    \end{enumerate}
    Furthermore, if $\ell_C(sw)<\ell_C(w)$ then there is a unique chamber $y\in X$ with $\delta(y,x)=s$ and $\delta^*(y,x_A)=sw$. Similarly, if $\ell_C(wt)<\ell_C(w)$ then there is a unique chamber $y_A\in X_A$ with $\delta_A(x_A,y_A)=t$ and $\delta^*(x,y_A)=wt$.
\end{lemma}

\begin{proof}
    The proof is essentially \cite[Lem.~5.139]{abramenko-brown08}. If $\ell_C(sw)>\ell_C(w)$ then $\delta^*(y,x_A)=sw$ by \ref{twisted:def:def:2}. So assume otherwise. By \ref{twisted:def:def:3} there exists $z\in X$ with $\delta(z,x)=s$ and $\delta^*(z,x_A)=sw$. If $y=z$ then (\ref{twisted:proj:lem:panels:1}) is satisfied, so assume otherwise. Since $\delta(z,x)=s$, and $\delta(y,x)=s$ by the hypothesis, one has $\delta(y,z)\in\{1,s\}$. But $y\neq z$ by assumption, so $\delta(y,z)=s$. Also, by assumption, one has $\ell_C(s(sw))>\ell_C(sw)$. Therefore $\delta^*(y,x_A)=s(sw)=w$ by \ref{twisted:def:def:2}. This proves (\ref{twisted:proj:lem:panels:1}). The argument for (\ref{twisted:proj:lem:panels:2}) is the same. The final statement is a repackaging of this argument.
\end{proof}

In short, if $x$ is not already the projection of $x_A$ onto the $s$-panel containing $x$, then the unique chamber $y$ of Lemma \ref{twisted:proj:lem:panels} is the projection. The analogous statement holds for the $t$-panel of $x_A$. Before extending the notion of projections to larger residues, the following analogue of \cite[Lem.~5.140]{abramenko-brown08} must be proved.

\begin{lemma}\label{twisted:proj:lem:galleries}
    Let $x,y\in X$ and $x_A,y_A\in X_A$. Suppose that
    \begin{equation*}
        u=\delta(x,y),\quad v=\delta^*(y,x_A),\quad\text{and}\quad w=\delta_A(x_A,y_A).
    \end{equation*}
    Then, for all $s_i\in S$ and $t_i\in S_A$, the following are true:
    \begin{enumerate}
        \item\label{twisted:proj:lem:galleries:1} If $u=s_1\dotsb s_l$ then $\delta^*(x,x_A)=s_{i_1}\dotsb s_{i_k}v$ for some subword $s_{i_1}\dotsb s_{i_k}$ of $s_1\dotsb s_l$.
        \item\label{twisted:proj:lem:galleries:2} If $w=t_1\dotsb t_l$ then $\delta^*(y,y_A)=vt_{i_1}\dotsb t_{i_k}$ for some subword $t_{i_1}\dotsb t_{i_k}$ of $t_1\dotsb t_l$.
        \item\label{twisted:proj:lem:galleries:3} If $\ell_C(uv)=\ell(u)+\ell_C(v)$ then $\delta^*(x,x_A)=uv$.
        \item\label{twisted:proj:lem:galleries:4} If $\ell_C(vw)=\ell_C(v)+\ell'(w)$ then $\delta^*(y,y_A)=vw$.
    \end{enumerate}
\end{lemma}

\begin{proof}
    For (\ref{twisted:proj:lem:galleries:1}) choose a gallery in $X$ passing from $x$ to $y$ described by the word $s_1\dotsb s_l$. Then repeated applications of Lemma \ref{twisted:proj:lem:panels}(\ref{twisted:proj:lem:panels:1}) gives the result. In a similar fashion, one deduces (\ref{twisted:proj:lem:galleries:2}) using Lemma \ref{twisted:proj:lem:panels}(\ref{twisted:proj:lem:panels:2}) and a gallery in $X_A$ passing from $x_A$ to $y_A$ described by the word $t_1\dotsb t_l$. For (\ref{twisted:proj:lem:galleries:3}) and (\ref{twisted:proj:lem:galleries:4}) consider the two cases for $C$ separately:
    \begin{itemize}
        \item Suppose that $C=A$. Let $u=s_1\dotsb s_d$ be a reduced decomposition in $W$. By Proposition \ref{cox:length:prop}(\ref{cox:length:prop:mult-simple}), multiplication by a simple reflection in $S$ changes $\ell_C$ by exactly one. Therefore, the hypothesis implies that 
        \begin{equation*}
        \ell_A(s_i\dotsb s_dv)>\ell_A(s_{i+1}\dotsb s_dv)\ \text{and}\ \ell_A(s_dv)>\ell_A(v)
        \end{equation*}
        for all $1\leq i\leq d-1$. The result now follows from repeated applications of \ref{twisted:def:def:2} on the gallery passing from $x$ to $y$ described by the word $s_1\dotsb s_d$.
        \item Suppose that $C=A^c$. Now, the hypothesis implies that 
        \begin{equation*}
        \ell_A(s_i\dotsb s_dv)<\ell_A(s_{i+1}\dotsb s_dv)\ \text{and}\ \ell_A(s_dv)<\ell_A(v)
        \end{equation*}
        for all $1\leq i\leq d-1$. This is equivalent to 
        \begin{equation*}
            \ell_{A^c}(s_i\dotsb s_dv)>\ell_{A^c}(s_{i+1}\dotsb s_dv)\ \text{and}\ \ell_{A^c}(s_dv)>\ell_{A^c}(v)
        \end{equation*}
        for all $1\leq i\leq d-1$. The same argument as before yields the result.
    \end{itemize}
    An analogous argument is made for (\ref{twisted:proj:lem:galleries:4}). Indeed, multiplication by a simple reflection in $S_A$ changes $\ell_C$ by exactly one using Proposition \ref{cox:refl-subgps:prop}. If $w=t_1\dotsb t_d$ is a reduced decomposition in $W_A$ then the argument follows identically.
\end{proof}

Let $J\subseteq S$ and $K\subseteq S_A$. Suppose that $R$ is a $J$-residue of $X$ and $R_A$ is a $K$-residue of $X_A$. Set
\begin{equation*}
    \delta^*(R,R_A)=\{\delta^*(x,x_A)\ |\ x\in R,\ x_A\in R_A\}.
\end{equation*}

\begin{lemma}\label{twisted:proj:lem:double-coset}
    Let $x\in R$ and $x_A\in R_A$. Then $\delta^*(R,R_A)=W_J\delta^*(x,x_A)(W_A)_K$.
\end{lemma}

\begin{proof}
    The proof follows \cite[Lem.~5.29]{abramenko-brown08}. Suppose that $\delta^*(x,x_A)=v$. If $y\in R$ and $y_A\in R_A$ then $\delta(y,x)\in W_J$ and $\delta_A(x_A,y_A)\in (W_A)_K$. Set $u=\delta(y,x)$ and $w=\delta_A(x_A,y_A)$. Applying Lemma \ref{twisted:proj:lem:galleries}(\ref{twisted:proj:lem:galleries:1}) and Lemma \ref{twisted:proj:lem:galleries}(\ref{twisted:proj:lem:galleries:2}) in succession implies that $\delta^*(y,y_A)=u'vw'$ for some subwords $u'$ of $u$ and $w'$ of $w$. Therefore $\delta^*(R,R_A)\subseteq W_J\delta^*(x,x_A)(W_A)_K$. This becomes an equality using \ref{twisted:def:def:3}.
\end{proof}

By Proposition \ref{ctd:prelim:cosets:prop:howlett}, if $W_J$ and $(W_A)_K$ are both finite then there are unique representatives 
\begin{equation*}
    w_A\in{}^JW_A^K\quad\text{and}\quad w_{A^c}\in{}^JW_{A^c}^K
\end{equation*}
of minimal twisted length $\ell_A(w_A)$ and $\ell_{A^c}(w_{A^c})$ in $\delta^*(R,R_A)$ respectively. In this case, the distances between chambers of $R$ and chambers of $R_A$ are bounded above and below, irrespective of the choice for $C$.

\begin{definition}
Let $R$ be a residue of $X$ and $R_A$ a residue of $X_A$. Then:
    \begin{itemize}
        \item The \emph{projection}, if it exists, of a chamber $x_A\in X_A$ onto $R$ is the unique chamber $\operatorname{proj}_R(x_A)$ of $R$ that minimises $d^C(-,x_A)$.
        \item The \emph{projection}, if it exists, of a chamber $x\in X$ onto $R_A$ is the unique chamber $\operatorname{proj}_{R_A}(x)$ of $R_A$ that minimises $d^C(x,-)$.
    \end{itemize}
\end{definition}

The projection of a chamber onto a residue may not always exist. This is already well-known for twin buildings, where projections across components exist only if the target residue is spherical. The crucial point here is that the distance from any chamber to a spherical residue is bounded below.

\begin{proposition}\label{twisted:proj:prop:residues}
    Let $x\in X$ and $x_A\in X_A$. Then:
    \begin{enumerate}
        \item\label{twisted:proj:prop:residues:1} If $d^C(R,x_A)$ is bounded below, then there exists a unique chamber $y\in R$ whose distance $d^C(y,x_A)$ to $x_A$ is minimal. Moreover, if $x\in R$ then
        \begin{equation*}
            \delta^*(x,x_A)=\delta(x,y)\,\delta^*(y,x_A)\qquad\text{and}\qquad d^C(x,x_A)=d(x,y)+d^C(y,x_A),
        \end{equation*}
        where $d$ denotes the gallery distance in $X$.
        \item\label{twisted:proj:prop:residues:2} If $d^C(x,R_A)$ is bounded below, then there exists a unique chamber $y_A\in R_A$ whose distance $d^C(x,y_A)$ to $x$ is minimal. Moreover, if $x_A\in R_A$ then
        \begin{equation*}
            \delta^*(x,x_A)=\delta^*(x,y_A)\,\delta_A(y_A,x_A)\qquad\text{and}\qquad d^C(x,x_A)=d^C(x,y_A)+d_A(y_A,x_A),
        \end{equation*}
        where $d_A$ denotes the gallery distance in $X_A$.
    \end{enumerate}
\end{proposition}

\begin{proof}
    The proof follows \cite[Lem.~5.149]{abramenko-brown08}. For (\ref{twisted:proj:prop:residues:1}), set $R_A=\{x_A\}$. Applying Lemma \ref{twisted:proj:lem:double-coset} to $R_A$ implies that $\delta^*(R,x_A)=W_J\delta^*(x,x_A)$ for any $x\in R$. By the hypothesis, there is a minimal element $v\in{}^JW_C$ in $\delta^*(R,x_A)$ which is unique by Proposition \ref{ctd:prelim:cosets:prop:howlett}. Choose $y\in R$ with $\delta^*(y,x_A)=v$. Fix $x\in R$ and set $u=\delta(x,y)\in W_J$. Applying Proposition \ref{cox:cosets:prop:posets}(\ref{cox:cosets:prop:posets:left}) to $u$ and $v$ gives 
    \begin{equation*}
        \ell_C(uv)=\ell(u)+\ell_C(v).
    \end{equation*}
    From Lemma \ref{twisted:proj:lem:galleries}(\ref{twisted:proj:lem:galleries:3}), one obtains
    \begin{equation*}
        \delta^*(x,x_A)=uv=\delta(x,y)\delta^*(y,x_A).
    \end{equation*}
    It follows that
    \begin{equation*}
        d^C(x,x_A)=\ell_C(uv)=\ell(u)+\ell_C(v)=d(x,y)+d^C(y,x_A).
    \end{equation*}
    Note that $\ell(u)>0$ provided $x\neq y$. Therefore, the chamber $y$ is the unique chamber of $R$ at minimal distance from $x_A$. This proves (\ref{twisted:proj:prop:residues:1}). The proof of (\ref{twisted:proj:prop:residues:2}) is analogous, using Proposition \ref{cox:cosets:prop:posets}(\ref{cox:cosets:prop:posets:right}) and Lemma \ref{twisted:proj:lem:galleries}(\ref{twisted:proj:lem:galleries:4}) instead.
\end{proof}

In general, a sufficient condition for the projection to exist is sphericity of the target residue. A more precise characterisation can be given in terms of initial sections. This rests on the following technical lemma. Write $T_J=T\cap W_J$ and $T_K=T_A\cap (W_A)_K$.

\begin{lemma}\label{twist:proj:lem:technical}
    Let $B$ be an initial section. Then
    \begin{equation*}
        \ell_B\ \text{bounded below on}\ W_J\iff |B\cap T_J|<\infty.
    \end{equation*}
    These conditions hold if and only if $B\cap T_J=N(u^{-1})$ for some $u\in W_J$. Similarly, one has
    \begin{equation*}
        \ell_C\ \text{bounded below on}\ (W_A)_K\iff |C\cap T_K|<\infty,
    \end{equation*}
    which occurs precisely when $C\cap T_K=N_A(u)$ for some $u\in (W_A)_K$.
\end{lemma}

\begin{proof}
    If $m=|B\cap T_J|<\infty$ then, for all $u\in W_J$, one has
    \begin{align*}
        \ell_B(u)&=\ell(u)-2|N(u^{-1})\cap B|,\\
        &=\ell(u)-2|N(u^{-1})\cap (B\cap T_J)|,\\
        &\geq\ell(u)-2m,\\
        &\geq-2m,
    \end{align*}
    using the fact that $N(u^{-1})\subseteq T_J$. Conversely, if $v\in W_J$ attains the minimum then $\ell_B(vt)>\ell_B(v)$ for all $t\in T_J$. Therefore $t\notin N(v^{-1})+B$ for all $t\in T_J$. It follows that
    \begin{equation*}
        (N(v^{-1})+B)\cap T_J=\emptyset.
    \end{equation*}
    This implies that
    \begin{equation*}
        B\cap T_J=N(v^{-1}).
    \end{equation*}
    using the fact that $N(v^{-1})\subseteq T_J$. The right-hand side of this equation is finite, proving the claim. The second statement follows from the same argument as above. Note that
    \begin{align*}
        \ell_C(u)&=\ell'(u)-2|N_A(u)\cap C|,\\
        &=\ell'(u)-2|N_A(u)\cap (C\cap T_K)|,
    \end{align*}
    for all $u\in (W_A)_K$, using Corollary \ref{cox:cosets:cor:lengths} and the fact that $N_A(u)\subseteq T_K$. For the converse, one obtains
    \begin{equation*}
        C\cap T_K=N_A(v)
    \end{equation*}
    for the minimal element $v\in (W_A)_K$.
\end{proof}

When $B\cap T_J=N(u^{-1})$ for some $u\in W_J$, one has
\begin{equation*}
    \ell_B(v)=\ell(vu^{-1})-\ell(u)
\end{equation*}
for every $v\in W_J$. A similar formula holds for $\ell_C$ and $u\in (W_A)_K$. In particular, the conditions of Lemma \ref{twist:proj:lem:technical} mean that, upon restriction to the appropriate parabolic subgroup, the twisted length function is a standard length function, possibly up to translation and a constant. 

\begin{proposition}\label{twisted:proj:prop:existence}
    Let $x\in X$ and $x_A\in X_A$. Suppose that $R$ is a $J$-residue of $X$ and $R_A$ is a $K$-residue of $X_A$. Then:
    \begin{equation*}
        \operatorname{proj}_{R_A}(x)\ \text{exists}\iff|C\cap T_K|<\infty.
    \end{equation*}
    Similarly, one has
    \begin{equation*}
        \operatorname{proj}_{R}(x_A)\ \text{exists}\iff|wC\cap T_J|<\infty
    \end{equation*}
    where $w=\delta^*(r,x_A)$ for any fixed $r\in R$.
\end{proposition}

\begin{proof}
    If $|C\cap T_K|$ is finite then $\ell_C$ is bounded below on $(W_A)_K$ by Lemma \ref{twist:proj:lem:technical}. Using Lemma 
    \ref{twisted:proj:lem:double-coset}, write $\delta^*(x,R_A)=w(W_A)_K$ for some $w\in W$. Let $w=s_1\dotsb s_d$ be a reduced decomposition in $W$. By Proposition \ref{cox:length:prop}(\ref{cox:length:prop:mult-simple}), multiplication by a simple reflection in $S$ changes $\ell_C$ by exactly one. Therefore, if $m=|C\cap T_K|<\infty$ then
    \begin{equation*}
        \ell_C(wv)\geq\ell_C(v)-d\geq-2m-\ell(w)
    \end{equation*}
    for all $v\in (W_A)_K$. This means that $\ell_C$ is bounded below on $w(W_A)_K$. Conversely, if $v\in w(W_A)_K$ minimises $\ell_C$ then 
    \begin{equation*}
        C\cap T_K=N(v^{-1})\cap T_K
    \end{equation*} 
    by the proof of Lemma \ref{twist:proj:lem:technical}. For the second equivalence, recall that
    \begin{equation*}
        \ell_C(uw)=\ell_C(w)+\ell_{wC}(u)
    \end{equation*}
    for all $u\in W_J$. Therefore, $\ell_C$ is bounded below on $W_Jw$ if and only if $\ell_{wC}$ is bounded below on $W_J$. The result follows from Lemma \ref{twist:proj:lem:technical} using $B=wC$.
\end{proof}

In ordinary buildings projections always exist because $C$ is empty. For twin buildings, where $C$ is the set of all reflections, projections exist precisely when $T_J$ (resp.~$T_K$) is finite. In other words, when the target is a spherical residue. Note that Proposition \ref{twisted:proj:prop:existence} does not exclude projections from $X$ onto the entirety of $X_A$ and vice versa. The latter case is less interesting. Indeed, if $x_A\in X_A$ then 
\begin{equation*}
    \operatorname{proj}_X(x_A)\ \text{exists}\iff C=N(w)\ \text{for some}\ w\in W
\end{equation*}
by Lemma \ref{twist:proj:lem:technical}. In particular, $C$ lies in the $W$-orbit of the empty set \cite[Lem.~2.11]{dyer-hecke93}. This means that $X$ and $X_A$ are isomorphic after possibly conjugating their types, cf.~Example \ref{twisted:def:ex-standard}(\ref{twisted:def:ex-standard3}). On the other hand, there are instances when $\operatorname{proj}_{X_A}$ exists and is non-trivial.

\begin{corollary}\label{twisted:proj:cor:finite}
    Let $x\in X$. If $|C\cap T_A|$ is finite then $\operatorname{proj}_{X_A}(x)$ exists. Moreover, for every residue $R_A$ of $X_A$, one has
    \begin{equation*}
        \operatorname{proj}_{R_A}(x)=\operatorname{proj}^{X_A}_{R_A}(\operatorname{proj}_{X_A}(x))
    \end{equation*}
    where $\operatorname{proj}^{X_A}_{R_A}(-)$ denotes the ordinary projection onto $R_A$ in $X_A$.
\end{corollary}

\begin{proof}
    The first statement follows immediately from Proposition \ref{twisted:proj:prop:existence}. For the second statement, set $p=\operatorname{proj}_{X_A}(x)$. By Proposition \ref{twisted:proj:prop:residues}(\ref{twisted:proj:prop:residues:2}), one has
    \begin{equation*}
        d^C(x,x_A)=d^C(x,p)+d_A(p,x_A)
    \end{equation*}
    for all $x_A\in X_A$. Set $p'=\operatorname{proj}_{R_A}^{X_A}(p)$. The ordinary gate property for buildings gives
    \begin{equation*}
        d_A(p,r_A)=d_A(p,p')+d_A(p',r_A)
    \end{equation*}
    for all $r_A\in R_A$. Substituting the latter equation into the former yields
    \begin{align*}
        d^C(x,r_A)&=d^C(x,p)+d_A(p,p')+d_A(p',r_A),\\
        &=d^C(x,p')+d_A(p',r_A),
    \end{align*}
    for all $r_A\in R_A$, where the second equality is an application of Proposition \ref{twisted:proj:prop:residues}(\ref{twisted:proj:prop:residues:2}) with $x_A=p'$. But $d_A(p',r_A)$ is non-negative, vanishing only when $r_A=p'$. Hence
    \begin{equation*}
        d^C(x,r_A)\geq d^C(x,p')
    \end{equation*} 
    for all $r_A\in R_A$, with equality precisely when $r_A=p'$. It follows that $p'=\operatorname{proj}_{R_A}(x)$, as required.
\end{proof}

In other words, if the hypothesis of Corollary \ref{twisted:proj:cor:finite} is satisfied then there exists a projection map from the building $X$ onto the satellite building $X_A$. This will be observed in the sequel for an important class of twisted buildings arising from non-standard positive systems of affine root systems, cf.~Example \ref{cox:initial-sections:ex}. However, it will not be true of all twisted buildings associated to positive systems in this way, as illustrated by the following example.

\begin{example}
    Let $(W,S)$ be the Coxeter system associated to the following diagram.
    \begin{center}
        \begin{dynkinDiagram}[
            Coxeter,
            label,
            root radius=0.1cm,
            edge length=1cm
            ]A5

        \node[above] at ($(root 1)!0.5!(root 2)$) {$\infty$};
        \node[above] at ($(root 4)!0.5!(root 5)$) {$\infty$};

        \end{dynkinDiagram}
    \end{center}
    Set $J=\{s_1,s_2\}$ and $K=\{s_4,s_5\}$. Then $W_J$ and $W_K$ are commuting infinite parabolic subgroups. Define the positive system
    \begin{equation*}
        \Psi=(\Delta^+\setminus\Delta^+_J)\cup(-\Delta^+_J)
    \end{equation*}
    where $\Delta_J^+$ is the standard positive system of the $J$-subsystem. Then
    \begin{equation*}
        A_\Psi=\{r_\alpha\ |\ \alpha\in\Delta^+\cap(-\Psi)\}=T_J.
    \end{equation*}
    The conjugation action of $J$ (resp.~$K$) stabilises (resp.~fixes) $T_J$. Hence $J\cup K\subseteq S_{A_\Psi}$. It follows that $T_J$ and $T_K$ are contained in $T_{A_\Psi}$. In other words, both
    \begin{equation*}
        A_\Psi\cap T_{A_\Psi}\quad\text{and}\quad T_{A_\Psi}+A_\Psi
    \end{equation*}
    are infinite, meaning a twisted building of type $(W,S,C)$ with $C\in\{A_\Psi,T+A_\Psi\}$ does not satisfy the hypothesis of  Corollary \ref{twisted:proj:cor:finite}.
\end{example}

The subsection concludes by defining and characterising projections between residues.

\begin{definition}
    Let $R$ be a residue of $X$ and $R_A$ a residue of $X_A$. Then:
    \begin{itemize}
        \item The \emph{projection} of $R_A$ onto $R$ is the set
        \begin{equation*}
            \operatorname{proj}_{R}R_A=\{\operatorname{proj}_{R}r_A\ |\ r_A\in R_A\}
        \end{equation*}
        whenever the appropriate chamber projections exist.
        \item The \emph{projection} of $R$ onto $R_A$ is the set
        \begin{equation*}
            \operatorname{proj}_{R_A}R=\{\operatorname{proj}_{R_A}r\ |\ r\in R\}
        \end{equation*}
        whenever the appropriate chamber projections exist.
    \end{itemize}
\end{definition}

The next result collects analogues of \cite[Lem.~5.36, Prop.~5.37, Ex.~5.168]{abramenko-brown08}.

\begin{proposition}\label{twisted:proj:prop:residue-residue}
    Let $R$ be a $J$-residue of $X$ and $R_A$ a $K$-residue of $X_A$. Suppose that $d^C$ is bounded below on $R\times R_A$. Write $w$ for the unique element of $\delta^*(R,R_A)$ minimising $\ell_C(w)$. Set
    \begin{equation*}
        L=J\cap wKw^{-1}\quad\text{and}\quad M=K\cap w^{-1}Jw.
    \end{equation*}
    Then:
    \begin{enumerate}
        \item\label{twisted:proj:prop:residue-residue:1} The projection of $R_A$ onto $R$, and vice versa, are given by
        \begin{align*}
            R'=\operatorname{proj}_{R}R_A&=\{x\in R\ |\ w\in\delta^*(x,R_A)\},\ \text{and}\\
            R_A'=\operatorname{proj}_{R_A}R&=\{x_A\in R_A\ |\ w\in\delta^*(R,x_A)\},
        \end{align*}
        respectively.
        \item\label{twisted:proj:prop:residue-residue:2} The sets $R'$ and $R_A'$ are residues of type $L$ and $M$, respectively.
        \item\label{twisted:proj:prop:residue-residue:3} The restrictions 
        \begin{equation*}
            f=\operatorname{proj}_{R_A}|_{R'}:R'\rightarrow R_A'\quad\text{and}\quad g=\operatorname{proj}_{R}|_{R_A'}:R_A'\rightarrow R'
        \end{equation*}
        are mutually inverse bijections. Moreover, $\delta^*(x,f(x))=w$ for all $x\in R'$.
        \item\label{twisted:proj:prop:residue-residue:4} For all $x,y\in R'$, one has
        \begin{equation*}
            \delta_A(f(x),f(y))=w^{-1}\delta(x,y)w.
        \end{equation*}
        In particular, $f$ is an isomorphism of buildings. 
    \end{enumerate}
\end{proposition}

\begin{proof}
    By the hypothesis, $\ell_C$ is bounded below on $\delta^*(R,R_A)$. If $w$ achieves this minimum then $w\in{}^JW_C^K$ and uniqueness follows from Proposition \ref{ctd:prelim:cosets:prop:howlett}. Suppose that $w=\delta^*(x,x_A)$ for some $x\in R$ and $x_A\in R_A$. The minimality and uniqueness of $w$ in Proposition \ref{twisted:proj:prop:residues} implies that $x=\operatorname{proj}_R(x_A)$ and $x_A=\operatorname{proj}_{R_A}(x)$. This proves the forward inclusions of (\ref{twisted:proj:prop:residue-residue:1}). For the reverse inclusion of the first equality, suppose that $x=\operatorname{proj}_R(x_A)$ for some $x_A\in R_A$. Set $v=\delta^*(x,x_A)$. By the definition of projection, $v\in{}^JW_C$. Using Proposition \ref{ctd:prelim:cosets:prop:howlett}, write $v=abwc$ such that $ab\in W_J$, $c\in (W_A)_K$, and $wc\in{}^JW_C$. Then $wc=(ab)^{-1}v\in W_Jv$ and $v,wc\in{}^JW_C$. By uniqueness, one must have $v=wc$ and thus $ab=1$. It follows that $v\in w(W_A)_K$ and so $w\in v(W_A)_K=\delta^*(x,R_A)$. The reverse inclusion of the second equality is proved in an analogous way, using the alternate decomposition $v=awb'c$ of Proposition \ref{ctd:prelim:cosets:prop:howlett}, where $a\in W_J$, $b'c\in(W_A)_K$, and $aw\in W_C^K$. For (\ref{twisted:proj:prop:residue-residue:2}), choose $y\in R'$ and $y_A\in R_A'$ with $\delta^*(y,y_A)=w$. Then $y=\operatorname{proj}_R(y_A)$ by (\ref{twisted:proj:prop:residue-residue:1}). Let $x\in R$ and set $u=\delta(x,y)$. By the gate property of Proposition \ref{twisted:proj:prop:residues}(\ref{twisted:proj:prop:residues:1}), one has
    \begin{equation*}
        \delta^*(x,y_A)=\delta(x,y)\delta^*(y,y_A)=uw.
    \end{equation*}
    Now
    \begin{equation*}
        \delta^*(x,R_A)=\delta^*(x,y_A)(W_A)_K=uw(W_A)_K
    \end{equation*}
    where the first equality follows from Lemma \ref{twisted:proj:lem:double-coset}. Therefore
    \begin{equation*}
        x\in R'\iff w\in uw(W_A)_K\iff u\in W_J\cap w(W_A)_Kw^{-1}.
    \end{equation*}
    By Proposition \ref{ctd:prelim:cosets:prop:kilmoyer}, this is equivalent to $u\in W_L$ where $L=J\cap wKw^{-1}$. In other words, $R'$ is the $L$-residue of $X$ containing $y$. The proof that $R_A'$ is a residue of type $M$ is analogous, using (\ref{twisted:proj:prop:residue-residue:1}), Proposition \ref{twisted:proj:prop:residues}(\ref{twisted:proj:prop:residues:2}), and the fact that
    \begin{equation*}
        (W_A)_K\cap w^{-1}W_Jw=(W_A)_{K\cap w^{-1}Jw},
    \end{equation*}
    respectively. For (\ref{twisted:proj:prop:residue-residue:3}), let $x\in R'$. Then, by (\ref{twisted:proj:prop:residue-residue:1}), there is $x_A\in R_A$ with $\delta^*(x,x_A)=w$. As at the start of the proof, this implies that $x=\operatorname{proj}_R(x_A)$ and $x_A=\operatorname{proj}_{R_A}(x)$. That is, $x_A\in R_A'$ and $g(f(x))=x$. Repeating the same argument with $y_A\in R_A'$ gives $f(g(x_A))=x_A$. Therefore, $f$ and $g$ are mutually inverse. Finally, suppose that $x,y\in R'$. Then
    \begin{equation*}
        \delta^*(x,f(y))=\delta(x,y)\delta^*(y,f(y))=\delta(x,y)w
    \end{equation*}
    by the first gate property of Proposition \ref{twisted:proj:prop:residues}. Applying the other gate property of Proposition \ref{twisted:proj:prop:residues} gives
    \begin{equation*}
        \delta^*(x,f(y))=\delta^*(x,f(x))\delta_A(f(x),f(y))=w\delta_A(f(x),f(y)).
    \end{equation*}
    It follows that
    \begin{equation*}
        \delta_A(f(x),f(y))=w^{-1}\delta(x,y)w,
    \end{equation*}
    which proves (\ref{twisted:proj:prop:residue-residue:4}).
\end{proof}

\subsection{Apartments}\label{twisted:apartments}
The existence of projection maps in twisted buildings motivates the following definition, inspired by the analogous notion in twin buildings.

\begin{definition}
    Let $Y\subseteq X$ and $Y_A\subseteq X_A$ be non-empty subsets. A pair $(Y, Y_A)$ is called \emph{convex} if $\operatorname{proj}_P(x)\in Y\cup Y_A$ for every chamber $x\in Y\cup Y_A$ and every panel $P\subseteq X\cup X_A$ that meets $Y\cup Y_A$.
\end{definition}

Just as for twin buildings, the definition can be simplified with an easy lemma.

\begin{lemma}\label{twisted:conv:lem}
    Let $Y\subseteq X$ and $Y_A\subseteq X_A$ be non-empty subsets. Then the pair $(Y, Y_A)$ is convex if and only if it is closed under taking projections onto residues meeting $Y\cup Y_A$ (whenever they exist).
\end{lemma}

\begin{proof}
    Let $x\in Y\cup Y_A$ be a chamber and $P\subseteq X\cup X_A$ a panel.
    If $x$ and $P$ lie in the same building then the result is well-known \cite[Lem.~5.45]{abramenko-brown08}. The remaining case is essentially \cite[Ex.~5.169]{abramenko-brown08}.
\end{proof}

If $(Y, Y_A)$ is convex then $Y$ is convex in $X$ and $Y_A$ is convex in $X_A$. This is immediate from the definition. The converse is not necessarily true and already fails for twin buildings. Another similarity is the fact that the intersection (taken component-wise) of a family of convex pairs is again convex, provided the intersections remain non-empty. 

\begin{definition}
    Let $(Y,Y_A,\epsilon^*)$ be a twisted building of type $(W,S,C)$. An \emph{isometry} of twisted buildings 
    \begin{equation*}
        (\phi,\phi_A):(X,X_A,\delta^*)\rightarrow(Y,Y_A,\epsilon^*)
    \end{equation*}
    is a pair of building isometries $\phi:X\rightarrow Y$ and $\phi_A:X_A\rightarrow Y_A$ satisfying 
    \begin{equation*}
        \epsilon^*(\phi(x),\phi_A(x_A))=\delta^*(x,x_A)
    \end{equation*}
    for all $x\in X$ and $x_A\in X_A$. If $\phi$ and $\phi_A$ are bijective then $(\phi,\phi_A)$ is an \emph{isomorphism}.
\end{definition}

A building isometry is necessarily injective, so an isometry of twisted buildings is always injective. The same is not true for surjectivity. In particular, an isometry $(\phi,\phi_A)$ is an isomorphism precisely when $\phi$ and $\phi_A$ are isomorphisms.

\begin{definition}
    A \emph{twisted apartment} of $(X,X_A,\delta^*)$ is a pair $(\Sigma,\Sigma_A)$ such that $\Sigma=\phi(W)$ and $\Sigma_A=\phi_A(W_A)$ for some isometry
    \begin{equation*}
        (\phi,\phi_A):(W,W_A,\delta^*_W)\rightarrow(X,X_A,\delta^*)
    \end{equation*}
    where $(W,W_A,\delta^*_W)$ is the thin twisted building of Example \ref{twisted:def:ex-thin}. 
\end{definition}

In other words, a twisted apartment $(\Sigma,\Sigma_A)$ consists of an apartment $\Sigma$ of $X$ and an apartment $\Sigma_A$ of $X_A$, subject to the compatibility condition
\begin{equation*}
    \delta^*(\phi(v),\phi_A(w))=v^{-1}w
\end{equation*}
for all $v\in W$ and $w\in W_A$.

\begin{lemma}\label{twisted:apartment:lem:convex}
    A twisted apartment $(\Sigma,\Sigma_A)$ is convex.
\end{lemma}

\begin{proof}
    Let $x$ be a chamber of $\Sigma\cup\Sigma_A$ and $P$ a panel that meets it. If $x$ and $P$ lie in the same building, then $\operatorname{proj}_P(x)\in\Sigma\cup\Sigma_A$ by convexity of building apartments \cite[Prop.~4.40]{abramenko-brown08}. Instead, assume that $x=\phi(u)$ for some $u\in W$ and suppose that $P\subseteq X_A$ is a $t$-panel through $\phi_A(w)$ meeting $\Sigma_A$, where $t\in S_A$ and $w\in W_A$. Then $\Sigma_A\cap P=\{\phi_A(w),\phi_A(wt)\}$. By Lemma \ref{twisted:proj:lem:panels}(\ref{twisted:proj:lem:panels:2}), this implies that $\delta^*(x,x_A)\in\{u^{-1}w,u^{-1}wt\}$ for all $x_A\in P$. But 
    \begin{equation*}
        |\ell_C(u^{-1}w)-\ell_C(u^{-1}wt)|=1
    \end{equation*}
    by Proposition \ref{cox:refl-subgps:prop} which means $\operatorname{proj}_P(x)\in \Sigma_A\cap P$ as required. The final case follows an analogous argument, using Lemma \ref{twisted:proj:lem:panels}(\ref{twisted:proj:lem:panels:1}) and Proposition \ref{cox:length:prop}(\ref{cox:length:prop:mult-simple}) instead.
\end{proof}

The remainder of the subsection is devoted to proving the common apartment property in twisted buildings. First, one must prove the following refinement of \cite[Thm.~5.73]{abramenko-brown08}.

\begin{lemma}\label{twisted:apartments:lem}
    Let $(X,\delta)$ be a building of type $(W,S)$. Suppose that there is a non-empty family of maps 
    \begin{equation*}
        f_i:X\rightarrow W
    \end{equation*}
    running over an index set $I$, such that, for all $x\in X$ and $s\in S$, the following properties hold:
    \begin{enumerate}
        \item If $y\in X$ with $\delta(x,y)=s$, then $f_i(y)\in\{f_i(x),f_i(x)s\}$.
        \item If $f_i(x)=w$ for every $i$, then there exists $y\in X$ with $\delta(x,y)=s$ and $f_i(y)=ws$ for every $i$.
    \end{enumerate}
    Then every isometry $\phi:V\rightarrow X$, with $V\subseteq W$ non-empty, satisfying
    \begin{equation*}
        f_i\circ\phi=\operatorname{id}_{V}
    \end{equation*}
    for all $i$, extends to a common isometric section $\widetilde{\phi}:W\rightarrow X$.
\end{lemma}

\begin{proof}
    Following \cite[Ex.~5.76]{abramenko-brown08}, consider the poset of pairs $(V',\phi')$ with $V\subseteq V'\subseteq W$ and $\phi':V'\rightarrow X$ an isometry extending $\phi$, satisfying $f_i\circ\phi'=\operatorname{id}_{V'}$, ordered by extension. Applying Zorn's lemma gives a maximal element $\widetilde{\phi}:D\rightarrow X$. By \cite[Thm.~5.73]{abramenko-brown08}, there is an isometry $\widehat{\phi}:W\rightarrow X$ extending $\widetilde{\phi}$. Observe that
    \begin{equation*}
        D=\{w\in W\ |\ f_i(\widehat{\phi}(w))=w\ \text{for all}\ i\in I\}.
    \end{equation*}
    The forward inclusion follows from the definition of $\widehat{\phi}$ and equality since $\widetilde{\phi}$ is maximal. In particular, $D$ is convex. Indeed, if $(w_0,\dotsb,w_m)$ is a minimal gallery of type $(s_1,\dotsb,s_m)$ with $w_0,w_m\in D$ then the first hypothesis implies that
    \begin{equation*}
        f_i(\widehat{\phi}(w_j))=w_0s_1^{\epsilon_1}\dotsb s_j^{\epsilon_j},\quad\epsilon_k\in\{0,1\}
    \end{equation*}
    for all $i\in I$ and $0\leq j\leq m$. If $j=m$ then $f_i(\widehat{\phi}(w_m))=w_m$ for all $i\in I$, by definition of $D$. By minimality, one must have $w_j\in D$ for all $j$. This proves the claim. Suppose for a contradiction that $D\neq W$. Since $D$ is non-empty and properly contained in $W$, there exists $v\in D$ and $s\in S$ such that $vs\notin D$. It follows that
    \begin{equation*}
        \ell(u^{-1}vs)=\ell(u^{-1}v)+1
    \end{equation*}
    for all $u\in D$. If this were not the case, then $vs$ would lie on a minimal gallery from $u$ to $v$, contradicting convexity. By the second hypothesis, there exists $y\in X$ with $\delta(\widetilde{\phi}(v),y)=s$ and $f_i(y)=vs$ for all $i\in I$. The building axioms in $X$ imply that
    \begin{equation*}
        \delta(\widetilde{\phi}(u),y)=u^{-1}vs
    \end{equation*}
    for all $u\in D$. In particular, setting $\widetilde{\phi}(vs)=y$ extends $\widetilde{\phi}$ to a common isometric section on $D\cup\{vs\}$, contradicting maximality. As such, $D=W$ and $\widetilde{\phi}$ is the desired common isometric section.
\end{proof}

The next result provides a strong extension property of isometries in twisted buildings.

\begin{theorem}\label{twisted:apartments:thm}
    Let $x_A\in X_A$. Suppose that $\phi:V\rightarrow X$ is an isometry from a non-empty subset $V\subseteq W$, satisfying 
    \begin{equation*}
        \delta^*(\phi(v),x_A)=v^{-1}
    \end{equation*}
    for all $v\in V$. Then there exist isometries
    \begin{equation*}
        \widetilde{\phi}:W\rightarrow X\quad\text{and}\quad\widetilde{\phi}_A:W_A\rightarrow X_A
    \end{equation*}
    such that
    \begin{equation*}
        \widetilde{\phi}|_V=\phi,\quad\widetilde{\phi}_A(1)=x_A,\quad\text{and}\quad\delta^*(\widetilde{\phi}(w),\widetilde{\phi}_A(w_A))=w^{-1}w_A
    \end{equation*}
    for all $w\in W$ and $w_A\in W_A$. Moreover, for a fixed $\widetilde{\phi}$, the map $\widetilde{\phi}_A$ is unique.
\end{theorem}

\begin{proof}
    To obtain the first isometry $\widetilde{\phi}$, apply Lemma \ref{twisted:apartments:lem} to $X$, $\phi$, and 
    \begin{equation*}
        f:X\rightarrow W,\ x\mapsto\delta^*(x,x_A)^{-1}.
    \end{equation*}
    The fact that $f$ satisfies the required hypotheses is immediate from Lemma \ref{twisted:proj:lem:panels}. By Lemma \ref{twisted:proj:lem:double-coset}, one has $\delta^*(\widetilde{\phi}(w),X_A)=w^{-1}W_A$. The isometry $\widetilde{\phi}_A$ is obtained from Lemma \ref{twisted:apartments:lem}, using $X_A$, $\phi_A:1\mapsto x_A$, and the maps
    \begin{equation*}
        f_w:X_A\rightarrow W_A,\ y_A\mapsto w\delta^*(\widetilde{\phi}(w),y_A)
    \end{equation*}
    for all $w\in W$. In particular, $f_w(\widetilde{\phi}_A(w_A))=w_A$ for all $w\in W$ and $w_A\in W_A$. Unravelling the definition, this gives
    \begin{equation*}
        \delta^*(\widetilde{\phi}(w),\widetilde{\phi}_A(w_A))=w^{-1}w_A
    \end{equation*}
    as required. It remains to show that the $f_w$ satisfy the hypotheses of Lemma \ref{twisted:apartments:lem}. The first hypothesis follows from Lemma \ref{twisted:proj:lem:panels}. The proof of the second hypothesis is more involved. Suppose that $y_A\in X_A$ with $f_w(y_A)=v$ for all $w\in W$. Fix $t\in S_A$. The goal is to find $z_A\in X_A$ with $\delta_A(z_A,y_A)=t$ and
    \begin{equation*}
        \delta^*(\widetilde{\phi}(w),z_A)=w^{-1}vt
    \end{equation*}
    for all $w\in W$. Then
    \begin{equation*}
        f_w(z_A)=w\delta^*(\widetilde{\phi}(w),z_A)=vt=f_w(y_A)t.
    \end{equation*}
    Set
    \begin{equation*}
        \alpha_\pm=\{w\in W\ |\ \ell_C(w^{-1}vt)-\ell_C(w^{-1}v)=\pm1\}.
    \end{equation*}
    In other words, one has
    \begin{equation*}
        w\in \alpha_+\iff t\notin N(v^{-1}w)+C.
    \end{equation*}
    It follows that $\alpha_+$ and $\alpha_-$ are roots of $W$ with common wall defined by $vtv^{-1}$. By Lemma \ref{twisted:proj:lem:panels}, if $z_A$ is $t$-adjacent to $y_A$ then
    \begin{equation*}
        \delta^*(\widetilde{\phi}(w),z_A)=w^{-1}vt
    \end{equation*}
    for all $w\in\alpha_+$. Pick $w_1\in\alpha_-$. By Lemma \ref{twisted:proj:lem:panels}, there exists a unique $t$-adjacent chamber $z_A$ of $y_A$ with
    \begin{equation*}
        \delta^*(\widetilde{\phi}(w_1),z_A)=w_1^{-1}vt.
    \end{equation*}
    This can be extended to all of the convex set $\alpha_-$ by induction on galleries. Assume that $w\in\alpha_-$ also satisfies this condition. Let $ws\in\alpha_-$ for some $s\in S$. By Lemma \ref{twisted:proj:lem:panels}, one has
    \begin{equation*}
        \delta^*(\widetilde{\phi}(ws),z_A)\in\{sw^{-1}v,sw^{-1}vt\}\cap\{w^{-1}vt,sw^{-1}vt\}.
    \end{equation*}
    This either takes value $sw^{-1}vt$ or
    \begin{equation*}
        sw^{-1}v=w^{-1}vt\iff vtv^{-1}=wsw^{-1}.
    \end{equation*}
    The latter equality implies that passage from $w$ to $ws$ crosses the common wall of $\alpha_+$ and $\alpha_-$ which is impossible. Hence
    \begin{equation*}
        \delta^*(\widetilde{\phi}(ws),z_A)=sw^{-1}vt=(ws)^{-1}vt
    \end{equation*}
    as required. For the final statement, observe that whenever $f_w(y_A)=v$ for all $w\in W$, the $t$-adjacent chamber $z_A$ satisfying $f_w(z_A)=vt$ is unique, as it is determined at $w_1$. But every choice of map $\widetilde{\phi}_A$ equals $x_A$ at $1$. As such, induction along galleries yields uniqueness on all of $W_A$.
\end{proof}

As a special case of Theorem \ref{twisted:apartments:thm}, one can obtain an analogue of the common apartment property known in buildings and twin buildings alike.

\begin{corollary}\label{twisted:apartments:cor:common}
    Let $x\in X$ and $x_A\in X_A$. Then $x$ and $x_A$ share a common twisted apartment.
\end{corollary}

\begin{proof}
    Consider the map
    \begin{equation*}
        \phi:\{w^{-1}\}\rightarrow X,\ w^{-1}\mapsto x,
    \end{equation*}
    where $w=\delta^*(x,x_A)$. Applying Theorem \ref{twisted:apartments:thm} to $\phi$ yields a twisted apartment $(\Sigma,\Sigma_A)$ where $\Sigma=\widetilde{\phi}(W)$ and $\Sigma_A=\widetilde{\phi}_A(W_A)$. By construction, one has
    \begin{equation*}
        x=\widetilde{\phi}(w^{-1})\in\Sigma\quad\text{and}\quad x_A=\widetilde{\phi}_A(1)\in\Sigma_A.
    \end{equation*}
\end{proof}

In contrast to twin buildings, two chambers in the same component of a twisted building need not share a common twisted apartment. This is illustrated in the next example.

\begin{example}\label{twisted:apartments:ex}
    Let $X$ be a thick affine tree and $X_A=\{x_A\}$ a point. Suppose that $\delta^*$ is a twisting map induced by a retraction from infinity, as in Example \ref{twisted:def:ex-retract}. By Lemma \ref{twisted:proj:lem:panels}, the map $\delta^*(-,x_A)$ takes two values on every panel. But $X$ is thick, so there must be distinct chambers $x,y\in X$ with $\delta^*(x,x_A)=\delta^*(y,x_A)$. Observe that $x$ and $y$ do not share a common twisted apartment. Indeed, any twisted apartment has the form $(\Sigma,\{x_A\})$ and is parameterised by a map $\phi:W\rightarrow \Sigma$ satisfying
    \begin{equation*}
        \delta^*(\phi(w),x_A)=w^{-1}
    \end{equation*}
    for all $w\in W$. In other words, $\delta^*(-,x_A)$ is injective on $\Sigma$. It follows that $x$ and $y$ cannot both be contained in $\Sigma$.
\end{example}

\begin{remark}\label{twisted:apartments:rmk}
    Let $(\Sigma,\Sigma_A)$ be a twisted apartment parameterised by $(\phi,\phi_A)$. Another difference from twin buildings is that a chamber $x\in\Sigma$ may have no chamber $x_A\in \Sigma_A$ for which $\delta^*(x,x_A)=1$. Indeed, if $x=\phi(w)$ then such a chamber exists if and only if $w\in W_A$. This has implications for any notion of opposition in twisted buildings which one may hope to define.
\end{remark}

\section{Affine twisted buildings from fa\c{c}ades}\label{aff}
This section constructs a class of twisted buildings associated to an affine building and its fa\c{c}ades. As a result, one obtains a different perspective on the chimney retractions of \cite{milicevic-gallery24}.

\subsection{Conventions}\label{aff:con}
Henceforth, let $(W,S)$ be an irreducible affine Coxeter system. In the sequel, it will be useful to pass to a model of an affine apartment as follows. Denote by $V$ the usual Euclidean affine space upon which $W$ acts. A \emph{wall} is a hyperplane in $V$ consisting of the fixed points for some reflection in $W$. Write $\mathcal{H}$ for the collection of all walls. An \emph{alcove} is the closure of a connected component in 
\begin{equation*}
    V\setminus\bigcup_{H\in\mathcal{H}}H.
\end{equation*}
Every alcove is a simplex and $W$ acts simply transitively on the set of alcoves. A face of an alcove is called a \emph{facet}. A facet of codimension one is called a \emph{panel}. 

The \emph{fundamental alcove} $c_f$ is the alcove bounded by the walls whose reflections lie in $S$. Fix a special vertex $v_0$ of $c_f$ as the origin of $V$. Let $W_0$ be its stabiliser in $W$ and $S_0$ its fixator in $S$. This defines a finite Coxeter system. Choose a reduced crystallographic root system $\Phi$ whose Weyl group is $W_0$ and affine Weyl group acts on $V$ as $W$. This is the convention taken in \cite[\S2.1]{milicevic-gallery24}. In this realisation, the real affine roots are
\begin{equation*}
    \Delta=\Phi+\mathbb{Z}\delta
\end{equation*}
where $\alpha+n\delta$ is the map on $V$ defined by $v\mapsto\langle\alpha,v\rangle+n$. Similarly, every wall of $V$ is of the form
\begin{equation*}
    H_{\alpha}=\{v\in V\ |\ \langle\alpha,v\rangle=0\}\quad\text{or}\quad H_{\alpha,n}=\{v\in V\ |\ \langle\alpha,v\rangle=n\}
\end{equation*}
for some root $\alpha\in\Phi$ and non-zero $n\in\mathbb{Z}$. Write $\Phi^+$ for the positive roots of $\Phi$ and $\Pi$ for the corresponding simple roots. 

\subsection{Fa\c{c}ades}\label{aff:fac}
Let $(X,\delta)$ be an affine building of type $(W,S)$. This subsection gives a brief account of a natural class of buildings at infinity associated to $X$. For more detail, see \cite{rousseau-book} and references therein. Throughout, $X$ is assumed to be \emph{chimney friendly}, meaning that any two chimney germs are contained in a common apartment.

Fix an apartment $\Sigma$ of $X$ and a vectorial model $|\Sigma|$ as in \S\ref{aff:con}. Let $\mathcal{H}_0$ be the collection of walls passing through the origin. A \emph{vectorial chamber} is the closure of a connected component in 
\begin{equation*}
    V\setminus\bigcup_{H\in\mathcal{H}_0}H.
\end{equation*}
The vectorial chambers in $V$ are non-empty simplicial cones, the faces of which are called \emph{vectorial facets}. A vectorial facet is determined, up to an element of $W_0$, by a subset $J$ in $\Pi$. More generally, a \emph{sector} is a cone of the form $x+C^v$ for some $x\in V$ and vectorial chamber $C^v$ in $V$. The faces of sectors are called \emph{sector faces}. These are cones of the form $\mathcal{F}=x+F^v$ for some vectorial facet $F^v$. In this instance, $x$ is called the \emph{base-point} and $F^v$ is called the \emph{direction}. More generally, one can define a \emph{chimney} to be the enclosure of the sum of a facet and a vectorial facet. In this case, the facet is called the \emph{basis} and the vectorial facet is again the \emph{direction}. A chimney is \emph{full} if it has a trivial pointwise stabiliser in $W$.

A pair of sectors in an apartment are \emph{parallel} if they contain a common subsector, i.e.~a sector of the same direction. This defines an equivalence relation, called \emph{parallelism}, on the collection of sectors. The collection $X^\infty$ of sectors whose base-points are special vertices, up to parallelism, is known as the \emph{spherical building at infinity}. Specifically, $X^\infty$ is a building of type $(W_0,S_0)$. It follows that every facet of $X^\infty$ can be represented by a sector face of $X$.

A similar relation can be defined for chimneys in the following way. The \emph{germ at infinity} $\operatorname{germ}_\infty(\mathcal{R})$ of a chimney $\mathcal{R}$ is the collection of chimneys $\mathcal{R}'$ such that $\mathcal{R}\cap\mathcal{R}'$ contains a chimney of the same direction, and generates the same affine subspace in any common apartment, as both $\mathcal{R}$ and $\mathcal{R}'$. For a pair of sector faces, it is enough to contain a subsector of the same direction. Let $F$ be a facet in $X^\infty$ of type $J\subseteq\Pi$. Now, the set
\begin{equation*}
    |X_F|=\{\operatorname{germ}_\infty(\mathcal{F})\ |\ \mathcal{F}\ \text{a sector face representing $F$}\}
\end{equation*}
defines an affine building called a \emph{fa\c{c}ade}, cf.~\cite[\S3.3.15]{rousseau-book}. Its chamber set $X_F$ is in one-to-one correspondence with the germs of full chimneys with direction $F$. The type of $|X_F|$ is isomorphic to the affine subsystem associated to $\Phi_J$. This will be denoted by $(W_F,S_F)$. Write
\begin{equation*}
    \delta_F:X_F\times X_F\rightarrow W_F
\end{equation*}
for the Weyl distance of $X_F$. 

\begin{remark}\label{aff:fac:rmk}
    The type of $X_F$ can be understood in terms of initial sections. In particular, it is the reflection subsystem in \S\ref{cox:refl-subgps} associated to the initial section
    \begin{equation*}
        A_F=\{r_{\alpha+n\delta}\ |\ \alpha\in\Phi^+\setminus\Phi_J,\ n\in\mathbb{Z}_{\geq0}\}.
    \end{equation*}
    To see this, note that $A_F$ is the initial section associated to the non-standard positive system
    \begin{equation*}
        \Psi_J=(-(\Phi^+\setminus\Phi_J)+\mathbb{Z}\delta)\cup\Phi_J^+\cup(\Phi_J+\mathbb{N}\delta).
    \end{equation*}
    Let $J_1,\dotsb,J_m$ be the connected components of $J$ and write $w_{0,i}$ for the highest root of $\Phi_{J_i}$. Then
    \begin{equation*}
        S_{A_F}=\{r_\alpha\ |\ \alpha\in J\}\cup\{r_{\delta-w_{0,i}}\ |\ 1\leq i\leq m\}
    \end{equation*}
    and therefore
    \begin{equation*}
        W_{A_F}=\langle r_{\alpha+n\delta}\ |\ \alpha\in\Phi_J,\ n\in\mathbb{Z}\rangle=W_F
    \end{equation*}
    as claimed. It will be important for later that $A_F\cap T_F=\emptyset$, cf.~Corollary \ref{twisted:proj:cor:finite}. Note that $\Psi_J$ differs from the representative of \cite{jakobsen-kac83} with the same notation.
\end{remark}

Let $z\in |X|$. Write $\mathcal{F}_z(F)$ for the sector face based at $z$ with ideal face $F$. In an apartment containing $z$ and $F$, one has $\mathcal{F}_z(F)=z+F^v$ where $F^v$ represents $F$. There is a natural \emph{(vectorial) projection} map
\begin{equation*}
    \pi_F:|X|\rightarrow |X_F|,\ z\mapsto\operatorname{germ}_\infty(\mathcal{F}_z(F))
\end{equation*}
sending a point to the germ of the sector face based there \cite[\S3.3.17]{rousseau-book}. For an apartment $\Sigma$ of $X$ with $F\subseteq\Sigma^\infty$, the map $\pi_F$ is the quotient by
\begin{equation*}
    E=\operatorname{span}(F^v)=\bigcap_{\alpha\in J}\ker\alpha.
\end{equation*}
Note that $W_F$ acts trivially on the factor of $E$. A point $z\in |X|$ in the interior of an alcove $x\in X$ is sent by $\pi_F$ to the interior of a unique chamber in the fa\c{c}ade. Denote this chamber by $\overline{x}\in X_F$. 

Let $F^-$ be a face of $X^\infty$ opposite to $F$. Write $\operatorname{Apt}(F,F^-)$ for the set of apartments $\Sigma$ in $X$ such that $F,F^-\subseteq\Sigma^\infty$. Then
\begin{equation*}
    |X_{F,F^-}|=\bigcup_{\Sigma\in \operatorname{Apt}(F,F^-)}|\Sigma|
\end{equation*}
is called an \emph{inner fa\c{c}ade} associated to $F$. This set admits two distinct affine building structures on the same system $\operatorname{Apt}(F,F^-)$ of apartments. The \emph{fine} structure is obtained by keeping the ambient walls in $X$. Therefore, its chamber set $X_{F,F^-}^\mathrm{f}$ is a subset of $X$. This defines a building of type $(W,S)$ that is not necessarily thick \cite[Thm.~3.3.17.1]{rousseau-book}. The corresponding Weyl distance is a restriction of $\delta_X$. By way of contrast, the \emph{coarse} structure forgets those walls which are not parallel to $F$. The result is another building structure on the inner fa\c{c}ade of type $(W_F,S_F)$, cf.~\cite[Thm.~3.3.14.1]{rousseau-book}. The associated chamber set $X_{F,F^-}^\mathrm{c}$ is indexed by alcoves of $X_F$. Specifically, every coarse chamber is a union of fine alcoves, called a \emph{double chimney}, which will be denoted by $D_{x_F}$ for some $x_F\in X_F$. In fact, $X_{F,F^-}^\mathrm{c}$ can be canonically identified with $X_F$. This follows from the isomorphism
\begin{equation*}
    |\Sigma|\simeq|\Sigma_F|\times E,
\end{equation*}
for all $\Sigma\in\operatorname{Apt}(F,F^-)$. In other words, the building $X_{F,F^-}^\mathrm{c}$ is \emph{non-essential} and the \emph{essentialisation} is given by the restriction
\begin{equation*}
    \pi_F:|X_{F,F^-}|\rightarrow |X_F|
\end{equation*}
which forgets the Euclidean factors parallel to $F$.

There is also a \emph{retraction} map 
\begin{equation*}
    \rho:|X|\rightarrow|X_{F,F^-}|
\end{equation*}
centred at $F$, which fixes the inner fa\c{c}ade pointwise and is type preserving for the fine structure. Furthermore, one has $\pi_F\circ\rho=\pi_F$. Details of the construction can be found in \cite[\S3.3.17]{rousseau-book}. Denote the induced map from $X$ to $X_{F,F^-}^\mathrm{f}$ by the same symbol. This will be central to the twisting map defined later. It is closely related to chimneys in the following way. Let $\operatorname{Apt}_F(x_F)$ denote the set of apartments of $X$ containing the full chimney germ represented by $x_F\in X_F$.

\begin{proposition}\label{aff:fac:prop:chimneys}
    Let $x\in X$ and $x_F\in X_F$. Then:
    \begin{enumerate}
        \item\label{aff:fac:prop:chimneys1} There exists an apartment $\Sigma\in\operatorname{Apt}_F(x_F)$ containing $x$.
        \item\label{aff:fac:prop:chimneys2} If $\Sigma\in\operatorname{Apt}_F(x_F)$ then there is a type-preserving map
        \begin{equation*}
            r_{\Sigma,x_F}:X\rightarrow \Sigma
        \end{equation*}
        called the \emph{chimney retraction}  centred at $x_F$. Moreover, if $\Sigma'\in\operatorname{Apt}_F(x_F)$ then
        \begin{equation*}
            \theta=r_{\Sigma',x_F}|_\Sigma:\Sigma\rightarrow\Sigma'
        \end{equation*} 
        is a type-preserving apartment isomorphism and 
        \begin{equation*}
            r_{\Sigma',x_F}=\theta\circ r_{\Sigma,x_F}.
        \end{equation*}
        \item\label{aff:fac:prop:chimneys3} If $\Sigma\in\operatorname{Apt}(F,F^-)$ contains $D_{x_F}$ then
        \begin{equation*}
            r_{\Sigma,x_F}=r_{\Sigma,D_{x_F}}\circ\rho
        \end{equation*}
        where $r_{\Sigma,D_{x_F}}$ is the apartment retraction of $X_{F,F^-}^\mathrm{c}$ centred at $D_{x_F}$. 
    \end{enumerate}
\end{proposition}

\begin{proof}
    (\ref{aff:fac:prop:chimneys1}) \cite[Prop.~3.15]{milicevic-gallery24}, (\ref{aff:fac:prop:chimneys2}) \cite[Cor.~3.17, Def.~3.18, Lem.~3.33]{milicevic-gallery24}, (\ref{aff:fac:prop:chimneys3}) \cite[\S3.3.18.4]{rousseau-book}.
\end{proof}

Note that the retraction $r_{\Sigma,D_{x_F}}$ coincides with any apartment retraction of $X_{F,F^-}^\mathrm{f}$ centred at an alcove contained in $D_{x_F}$. This will be used in the sequel without further reference.

\subsection{A twisting map}\label{aff:twist}
As before, let $(X,\delta)$ be an affine building of type $(W,S)$. Fix a facet $F\subseteq X^\infty$ of type $J$. By Remark \ref{aff:fac:rmk}, the type of $X_F$ is the reflection subsystem associated to $A_F$. For brevity, write $\ell_F=\ell_{A_F}$ and $\ell'_F=\ell'_{A_F}$. The goal of this subsection is to define a twisting map between the building $X$ and the fa\c{c}ade $X_F$. For the most part, this will be determined by the retraction
\begin{equation*}
    \rho:X\rightarrow X_{F,F^-}^\mathrm{f}\subseteq X.
\end{equation*}  
However, one cannot directly compare alcoves in $X_{F,F^-}^\mathrm{f}$ with alcoves in $X_F$. Instead, identify $X_F$ with $X_{F,F^-}^\mathrm{c}$ as in \S\ref{aff:fac}. Keeping the notation as before, it remains to choose an alcove of $X$ in every $D_y$ in a well-defined manner. This can be achieved using the \emph{folding retraction}
\begin{equation*}
    \theta_D:|X_{F,F^-}^\mathrm{c}|\rightarrow D
\end{equation*}
for a fixed alcove $D$. Roughly speaking, this is the retraction centred at $D$ followed by a Weyl-translation in the resulting apartment, cf.~\cite[\S2.1.7]{rousseau-book}. Importantly, the restriction
\begin{equation*}
    \theta_D|_{D'}:D'\rightarrow D
\end{equation*}
is an isometry which respects the underlying fine structure. In general, $\theta_D$ is invariant under Weyl isomorphisms of the coarse structure, i.e.~apartment retractions.

Fix an apartment $\Sigma_0$ of $X_{F,F^-}$ and identify it with the model apartment of \S\ref{aff:con} such that $F$ corresponds to the negative standard vectorial facet of type $J$. Under this identification, set $x'=c_f$, $x_F'=\overline{c_f}$, and $D=D_{\overline{c_f}}$. 

\begin{lemma}\label{aff:twist:lem:j}
    Keeping the choices above, the map
    \begin{equation*}
        j:X_F\rightarrow X_{F,F^-}^\mathrm{f},\ x_F\mapsto\theta_D|_{D_{x_F}}^{-1}(x')
    \end{equation*}
    is well-defined and satisfies
    \begin{equation*}
        \overline{j(x_F)}=x_F\quad\text{and}\quad\delta(j(x_F),j(y_F))=\delta_F(x_F,y_F)
    \end{equation*}
    for all $x_F,y_F\in X_F$. Moreover, if $\Sigma$ is an apartment of $X_{F,F^-}$ containing $D_{x_F}$ then
    \begin{equation*}
        r_{\Sigma,D_{x_F}}\circ j=j\circ r_{\Sigma_F,x_F}
    \end{equation*}
    where $r_{\Sigma_F,{x_F}}$ is the retraction of $X_F$ onto $\Sigma_F$ centred at $x_F$.
\end{lemma}

\begin{proof}
    By construction, $x'$ is an alcove with respect to the fine structure. Let $x_F\in X_F$. Then $j(x_F)\subseteq D_{x_F}$ by the discussion above. The fact that $\overline{j(x_F)}=x_F$ is immediate, since $\pi_F$ sends the interior of $D_{x_F}$ onto $x_F$. Let $y_F\in X_F$. Suppose that $\Sigma$ is an apartment of $X_{F,F^-}$ containing $D_{x_F}$ and $D_{y_F}$. Choose an apartment $\Sigma'$ containing $D$ and $D_{x_F}$ along with isomorphisms
    \begin{equation*}
        \alpha:\Sigma\rightarrow\Sigma'\quad\text{and}\quad\beta:\Sigma'\rightarrow\Sigma_0
    \end{equation*}
    fixing $D_{x_F}$ and $D$ pointwise, respectively \cite[Thm.~2.1.8(d)]{rousseau-book}. Set $u=\delta_F(x_F',x_F)$. Then
    \begin{equation*}
        \varphi=u^{-1}\circ\beta\circ\alpha:\Sigma\rightarrow\Sigma_0
    \end{equation*}
    is an apartment isomorphism in the coarse structure extending $\theta_D|_{D_{x_F}}$. It is straightforward to see that $\varphi(j(x_F))=x'$. Set $v=\delta_F(x_F,y_F)$. By construction, one has $\varphi(D_{y_F})=vD$ and
    \begin{equation*}
        \theta_D(\varphi(j(y_F)))=\theta_D(j(y_F))=x'
    \end{equation*}
    where the first equality follows from the fact that $\theta_D$ is invariant under coarse Weyl isomorphisms. Note that $\varphi(j(y_F))\subseteq vD$ since $j(y_F)\subseteq D_{y_F}$. But on $vD$ the map $\theta_D$ is just multiplication by $v^{-1}$. Therefore, one has $\varphi(j(y_F))=vx'$. Observe that $\varphi$ preserves fine types, which implies that
    \begin{equation*}
        \delta(j(x_F),j(y_F))=\delta(c_f,vc_f)=v=\delta_F(x_F,y_F)
    \end{equation*}
    as required, where the first equality uses $x'=c_f$. For the final statement, set $z_F=r_{\Sigma_F,x_F}(y_F)$. The restriction of $r_{\Sigma,D_{x_F}}$ to a common apartment of $D_{x_F}$ and $D_{y_F}$ is an isomorphism. In this context, $D_{y_F}$ is sent to $D_{z_F}$ and the image of $j(y_F)$ is a fine alcove contained in $D_{z_F}$. By the same argument as before, one has
    \begin{equation*}
        \theta_D(r_{\Sigma,D_{x_F}}(j(y_F)))=\theta_D(j(y_F))=x'.
    \end{equation*}
    But $j(z_F)$ is the unique fine alcove in $D_{z_F}$ sent to $x'$ under $\theta_D$. Hence
    \begin{equation*}
        r_{\Sigma,D_{x_F}}(j(y_F))=j(z_F)=j(r_{\Sigma_F,x_F}(y_F)).
    \end{equation*}
\end{proof}

\begin{definition}\label{aff:twist:def}
    For $x\in X$ and $x_F\in X_F$, define
    \begin{equation*}
        \delta^*(x,x_F)=\delta(\rho(x),j(x_F))\quad\text{and}\quad\delta^*(x_F,x)=\delta^*(x,x_F)^{-1}.
    \end{equation*}
\end{definition}

Thus, one has a $W$-valued map
\begin{equation*}
    \delta^*:(X\times X_F)\cup(X_F\times X)\rightarrow W
\end{equation*}
where the union of chamber sets is taken to be disjoint. The map $\delta^*$ can be understood via chimney retractions as follows. Let $x\in X$ and $x_F\in X_F$. Suppose that $\Sigma$ is an apartment of $X_{F,F^-}$ containing $\rho(x)$ and $D_{x_F}$. By Proposition \ref{aff:fac:prop:chimneys}, one has
\begin{equation*}
    \delta^*(x,x_F)=\delta(r_{\Sigma,x_F}(x),j(x_F)).
\end{equation*}
In particular, $\delta^*(-,x_F)$ is a (normalised) chimney retraction centred at $x_F$. 

\begin{remark}
    The fixed model preceding Lemma \ref{aff:twist:lem:j} follows the standard $J$-chimney conventions of \cite{milicevic-gallery24}. A different choice of fine alcove $uc_f$ in $D$ would replace $\delta^*(x,x_F)$ by $\delta^*(x,x_F)u$ for all $x\in X$ and $x_F\in X_F$. In order to maintain compatibility in Lemma \ref{aff:twist:lem:j}, the corresponding fa\c{c}ade distance should be $u^{-1}\delta_Fu$ and then, jumping forward, the resulting twisting map will have type $u^{-1}A_F$.
\end{remark}

Fix an apartment $\Sigma$ of $X_{F,F^-}$ containing $D_{x_F}$. The chimney represented by $x_F$ can be understood as the collection of half-apartments of $\Sigma$ which contain a representative of its germ, cf.~\cite[Rmk.~3.2]{milicevic-gallery24}. It defines a \emph{chimney orientation} $\phi_{\Sigma,x_F}$ which associates a sign to every alcove $x\in\Sigma$ and panel $p$ contained in it. Specifically, one has $\phi_{\Sigma,x_F}(x,p)=+$ if and only if the half-apartment of $\Sigma$ bounded by $\operatorname{supp}(p)$ and containing $x$ does not belong to this collection, cf.~\cite[Def.~3.24]{milicevic-gallery24}.

\begin{lemma}\label{aff:twist:prop:orientation}
    Let $x\in X$ and $x_F\in X_F$ with $w=\delta(x,j(x_F))$. Suppose that $x\in \Sigma$ for some apartment $\Sigma$ of $X_{F,F^-}$ containing $D_{x_F}$. Then
    \begin{equation*}
        \phi_{\Sigma,x_F}(x,p)=+\iff s\in wA_F\iff \ell_F(sw)<\ell_F(w)
    \end{equation*}
    for any $s$-panel $p$ of $x$, where $s\in S$. 
\end{lemma}

\begin{proof}
    Let $f:V\rightarrow|\Sigma|$ be the type-preserving
    apartment isomorphism sending $c_f$ to $j(x_F)$ and the standard $J$-chimney to the chimney represented by $x_F$. The half-apartments of the former chimney are given by
    \begin{equation*}
        V_\alpha=\{v\in V\ |\ \alpha(v)\geq0\}
    \end{equation*}
    for some $\alpha\in\Psi_J$, those of the latter chimney by $f(V_\alpha)$. Observe that
    \begin{equation*}
        x=f(w^{-1}c_f)\quad\text{and}\quad p=f(w^{-1}p_s)
    \end{equation*}
    where $p_s$ is the $s$-panel of $c_f$. Therefore, the half-apartment bounded by $\operatorname{supp}(p)$ and containing $x$ is $f(V_{w^{-1}\alpha_s})$ where $\alpha_s$ is the simple root associated to $s$. It follows that
    \begin{equation*}
        \phi_{\Sigma,x_F}(x,p)=+\iff w^{-1}\alpha_s\in-\Psi_J\iff s\in A_{w\Psi_J}=wA_F
    \end{equation*} 
    where the second equivalence is by definition, cf.~Example \ref{cox:initial-sections:ex}. The final equivalence in terms of $\ell_F$ is Proposition \ref{cox:length:prop}(\ref{cox:length:prop:mult-char}).
\end{proof}

As a consequence of Lemma \ref{aff:twist:prop:orientation} one obtains a left-hand rule for $\delta^*$ in terms of $\ell_F$.

\begin{proposition}\label{aff:twist:prop:left}
    Let $x\in X$ and $x_F\in X_F$. Suppose that $\delta^*(x,x_F)=w$. If $y\in X$ with $\delta(y,x)=s$ for some $s\in S$, then $\delta^*(y,x_F)\in\{w,sw\}$. Furthermore, if $\ell_F(sw)>\ell_F(w)$ then $\delta^*(y,x_F)=sw$.
\end{proposition}

\begin{proof}
    Let $\Sigma$ be an apartment of $X_{F,F^-}$ containing $D_{x_F}$. Note that $r_{\Sigma,D_{x_F}}$ coincides with the fine chamber retraction of $X_{F,F^-}$ centred at $j(x_F)\subseteq D_{x_F}$. In particular, it preserves the Weyl distance $\delta(-,j(x_F))$. It follows that
    \begin{align*}
        \delta^*(z,x_F)&=\delta(\rho(z),j(x_F)),\\
        &=\delta((r_{\Sigma,D_{x_F}}\circ\rho)(z),j(x_F)),\\
        &=\delta(r_{\Sigma,x_F}(z),j(x_F)),
    \end{align*}
    for all $z\in X$, where the final equality is Proposition \ref{aff:fac:prop:chimneys}(\ref{aff:fac:prop:chimneys3}). Since $r_{\Sigma,x_F}$ is type-preserving and $\Sigma$ is thin, one deduces that $\delta^*(y,x_F)\in\{w,sw\}$. When $\delta^*(y,x_F)=w$ the images of $x$ and $y$ under $r_{\Sigma,x_F}$ are equal. Let $p$ be the common panel of $x$ and $y$. Then 
    \begin{equation*}
        \phi_{\Sigma,x_F}(r_{\Sigma,x_F}(x),r_{\Sigma,x_F}(p))=\phi_{\Sigma,x_F}(r_{\Sigma,x_F}(y),r_{\Sigma,x_F}(p))=+
    \end{equation*}
    by the proof of \cite[Prop.~3.32]{milicevic-gallery24}. The final statement now follows from Lemma \ref{aff:twist:prop:orientation}.
\end{proof}

It remains to prove the corresponding right-hand rule. This will be made possible by the next lemma, which should be seen as a gate property. Recall from \S\ref{aff:fac} that $x\in X$ is projected under $\pi_F$ to a unique chamber $\overline{x}$ in $X_F$. Moreover, one has $\pi_F\circ\rho=\pi_F$.

\begin{lemma}\label{aff:twist:lem:vec-gate}
    Let $x\in X$ and $x_F\in X_F$. Then
    \begin{equation*}
        \delta^*(x,x_F)=\delta^*(x,\overline{x})\delta_F(\overline{x},x_F).
    \end{equation*}
    Moreover, $\delta^*(x,\overline{x})$ is of minimal twisted length $\ell_F$ in its left coset of $W_F$. As a result, one has
    \begin{equation*}
        \ell_F(\delta^*(x,x_F))=\ell_F(\delta^*(x,\overline{x}))+\ell_F'(\delta_F(\overline{x},x_F)).
    \end{equation*}
\end{lemma}

\begin{proof}
    By the discussion above, $\rho(x)$ and $x$ project to the same fa\c{c}ade chamber $\overline{x}$. Hence $\rho(x)$ is contained in $D_{\overline{x}}$. By definition, $j(\overline{x})$ is also contained in $D_{\overline{x}}$. Let $\Sigma$ be an apartment of $X_{F,F^-}$ containing $D_{\overline{x}}$ and $D_{x_F}$. Then $\Sigma$ contains $\rho(x)$, $j(\overline{x})$, and $j(x_F)$. Since $\Sigma$ is thin, it follows that
    \begin{equation*}
        \delta(\rho(x),j(x_F))=\delta(\rho(x),j(\overline{x}))\delta(j(\overline{x}),j(x_F)).
    \end{equation*}
    This can be rewritten as
    \begin{equation*}
        \delta^*(x,x_F)=\delta^*(x,\overline{x})\delta_F(\overline{x},x_F)
    \end{equation*}
    using the definition of $\delta^*$ and Lemma \ref{aff:twist:lem:j}. As in the proof of Lemma \ref{aff:twist:lem:j}, there is a coarse apartment isomorphism 
    \begin{equation*}
        \varphi:\Sigma\rightarrow\Sigma_0
    \end{equation*}
    which sends $D_{\overline{x}}$ to $D$ and $j(\overline{x})$ to $x'$, preserving fine types. Write $\varphi(\rho(x))=wx'$ for some $w\in W$. Then
    \begin{align*}
        \delta^*(x,\overline{x})&=\delta(\rho(x),j(\overline{x})),\\
        &=\delta(wx',x'),\\
        &=w^{-1}.
    \end{align*}
    Note that $wx'$ and $x'$ both lie in $D$. Therefore, under the fixed identification of $\Sigma_0$ with $V$ such that $x'=c_f$, no wall of $W_F$ separates them. Equivalently, $N(w)\cap T_F=\emptyset$. But $A_F\cap T_F$ is also empty, meaning
    \begin{equation*}
        (N(w)+A_F)\cap T_F=\emptyset.
    \end{equation*}
    This implies the minimality of $w^{-1}=\delta^*(x,\overline{x})$ required, cf.~\S\ref{cox:cosets}. The final length formula is an application of Proposition \ref{cox:cosets:prop:posets}(\ref{cox:cosets:prop:posets:right}).
\end{proof}

The similarity of Lemma \ref{aff:twist:lem:vec-gate} and Corollary \ref{twisted:proj:cor:finite} is not a coincidence. After realising $\delta^*$ as a twisting map, it will be proved that $\pi_F$ induces a projection map of twisted buildings onto the satellite building $X_F$. This will be returned to after completing the aim of this section.

\begin{theorem}\label{aff:twist:thm}
    The triple $(X,X_F,\delta^*)$ is a twisted building of type $(W,S,A_F)$.
\end{theorem}

\begin{proof}
    Let $x\in X$ and $x_F\in X_F$ with $w=\delta^*(x,x_F)$. The axiom \ref{twisted:def:def:1} holds by definition. By Proposition \ref{aff:twist:prop:left}, the axiom \ref{twisted:def:def:2a} holds also. Using Lemma \ref{aff:twist:lem:vec-gate}, write $w=uv$ where $u=\delta^*(x,\overline{x})$ and $v=\delta_F(\overline{x},x_F)$. If $\delta_F(x_F,y_F)=t\in S_F$ for some $y_F\in X_F$ then
    \begin{equation*}
        \delta_F(\overline{x},y_F)\in\{v,vt\}
    \end{equation*}
    with $\delta_F(\overline{x},y_F)=vt$ whenever $\ell_F'(vt)>\ell_F'(v)$. Also by Lemma \ref{aff:twist:lem:vec-gate}, one has
    \begin{equation*}
        \ell_F(wt)-\ell_F(w)=\ell_F'(vt)-\ell_F'(v).
    \end{equation*}
    Therefore, $\ell_F(wt)>\ell_F(w)$ implies that $\delta^*(x,y_F)=uvt=wt$, proving axiom \ref{twisted:def:def:2b}. For \ref{twisted:def:def:3a}, fix an apartment $\Sigma$ of $X_{F,F^-}$ containing $D_{x_F}$. Let $\Sigma'$ be an apartment of $X$ containing $x$ and the chimney germ represented by $x_F$. Then 
    \begin{equation*}
        r_{\Sigma,x_F}|_{\Sigma'}:\Sigma'\rightarrow\Sigma
    \end{equation*}
    is a type-preserving isomorphism. Suppose that $s\in S$. If $y\neq x$ is the unique chamber of $\Sigma'$ containing the $s$-panel of $x$, then $\delta(y,x)=s$ and 
    \begin{align*}
        \delta^*(y,x_F)&=\delta(\rho(y),j(x_F)),\\
        &=\delta(r_{\Sigma,x_F}(y),j(x_F)),\\
        &=s\delta(r_{\Sigma,x_F}(x),j(x_F)),\\
        &=sw,
    \end{align*}
    by the discussion immediately after Definition \ref{aff:twist:def}. For all $t\in S_F$, the existence axiom in $X_F$ supplies a chamber $y_F$ with $\delta_F(x_F,y_F)=t$ and $\delta_F(\overline{x},y_F)=vt$. Applying Lemma \ref{aff:twist:lem:vec-gate} to $x$ and $y_F$ yields $\delta^*(x,y_F)=wt$. This proves \ref{twisted:def:def:3b} and the result follows.
\end{proof}

By Theorem \ref{aff:twist:thm}, there are well-defined notions of projections between residues in an affine building and its fa\c{c}ades. The next corollary shows that the projection from the affine building to its fa\c{c}ade is a combinatorial analogue of the vectorial projection in \S\ref{aff:fac}.

\begin{corollary}\label{aff:twist:cor:proj}
    The vectorial projection $\pi_F$ agrees with the twisted building projection map on chambers. That is, for all $x\in X$, one has
    \begin{equation*}
        \operatorname{proj}_{X_F}(x)=\overline{x}.
    \end{equation*} 
    Moreover, for every residue $R_F\subseteq X_F$, one has
    \begin{equation*}
        \operatorname{proj}_{R_F}(x)=\operatorname{proj}_{R_F}^{X_F}(\overline{x})
    \end{equation*}
    where $\operatorname{proj}^{X_F}_{R_F}(-)$ is the ordinary projection onto $R_F$ in $X_F$.
\end{corollary}

\begin{proof}
    The map $\operatorname{proj}_{X_F}(-)$ is well-defined by Corollary \ref{twisted:proj:cor:finite} since $A_F\cap T_{F}=\emptyset$. By Lemma \ref{aff:twist:lem:vec-gate} and the uniqueness of projections, $\overline{x}$ must coincide with $\operatorname{proj}_{X_F}(x)$. The remaining statement is an application of Corollary \ref{twisted:proj:cor:finite}.
\end{proof}

\appendix
\section{Proofs of Dyer's results}\label{ap:proofs}
This appendix records the proofs of \S\ref{cox}. As mentioned before, these proofs were communicated to the author by Dyer privately, and are reproduced here with permission.

Let $(W,S)$ be a Coxeter system and $A$ an initial section. Before proving Proposition \ref{cox:cosets:prop:posets}, one proves the following intermediate results.

\begin{lemma}\label{cox:cosets:lem:poset-isom-right}
    Let $J\subseteq S$. Suppose that $w=uv$ where $u\in W_J$ and $v\in W$. Then
    \begin{equation*}
        \ell_A(w)=\ell_{vA\cap W_J}(u)+\ell_A(v).
    \end{equation*}
    In particular, the map $W_J\rightarrow W_Jv$ defined by $u\mapsto uv$ is an isomorphism of posets, where $W_J$ has the order $\leq_{vA\cap W_J}$ and $W_Jv$ has the order $\leq_A$ induced from $W$.
\end{lemma}
\begin{proof}
    Let $t$ be a reflection of $W_J$. By Proposition \ref{cox:length:prop}(\ref{cox:length:prop:mult-char}), one has
    \begin{equation*}
        \ell_A(tw)<\ell_A(w)\iff t\in wA.
    \end{equation*}
    for all $w\in W$. But
    \begin{equation*}
        wA=u(vA)=N(u)+u(vA)u^{-1},
    \end{equation*}
    from which it follows that
    \begin{align*}
        \ell_A(tw)<\ell_A(w)&\iff t\in N(u)+u(vA)u^{-1},\\
        &\iff u^{-1}tu\in u^{-1}N(u)u+vA,\\
        &\iff u^{-1}tu\in N(u^{-1})+vA.
    \end{align*}
     Observe that $u^{-1}tu\in W_J$, so the latter condition can be rewritten as
     \begin{equation*}
        u^{-1}tu\in (N(u^{-1})\cap W_J)+(vA\cap W_J).
     \end{equation*}
     Altogether, one has
     \begin{align*}
        \ell_A(tw)<\ell_A(w)&\iff u^{-1}tu\in (N(u^{-1})\cap W_J)+(vA\cap W_J).\\
        &\iff t\in (N(u)\cap W_J)+u(vA\cap W_J)u^{-1}.
    \end{align*}
    Since $W_J$ is a standard parabolic subgroup, the standard length function on $W_J$ coincides with the standard length function on $W$. In particular, $N(u)\cap W_J$ is the inversion set of $u$ in $W_J$. Therefore
    \begin{equation*}
        (N(u)\cap W_J)+u(vA\cap W_J)u^{-1}=u(vA\cap W_J)
    \end{equation*}
    and thus $\ell_A(tw)<\ell_A(w)$ if and only if $\ell_{vA\cap W_J}(tu)<\ell_{vA\cap W_J}(u)$ by Proposition \ref{cox:length:prop}(\ref{cox:length:prop:mult-char}). The length additivity equality follows from induction on $\ell(u)$. For the base case, if $\ell(u)=0$ then $u$ is the identity in $W$ and the equality holds by Proposition \ref{cox:length:prop}(\ref{cox:length:prop:identity}). Assume the equality is true for $\ell(u)=n$. Suppose $t\in J$ with $\ell(tu)=\ell(u)+1$. By Proposition \ref{cox:length:prop}(\ref{cox:length:prop:mult-simple}), if $t\in J\subseteq S$ then $\ell_{vA\cap W_J}(tu)=\ell_{vA\cap W_J}(u)\pm1$, from which it follows that $\ell_A(tw)=\ell_A(w)\pm1$ (with compatible same signs) by the discussion above. Therefore
    \begin{align*}
        \ell_A(tw)&=\ell_A(w)\pm1,\\
        &=\ell_{vA\cap W_J}(u)+\ell_A(v)\pm1,\\
        &=\ell_{vA\cap W_J}(tu)+\ell_A(v),
    \end{align*}
    where the induction hypothesis has been used in the second equality. This proves the length additivity equation. The second claim now follows.
\end{proof}
\begin{lemma}\label{cox:cosets:lem:poset-isom-left}
    Suppose that $w=uv$ where $u^{-1}\in X_{W_A}$ and $v\in W_A$. Then
    \begin{equation*}
        \ell_A(w)=\ell_A(u)+\ell_A'(v^{-1}).
    \end{equation*}
    In particular, the map $W_A\rightarrow wW_A$ defined by $v^{-1}\mapsto uv$ is an isomorphism of posets, where $W_A$ has the order $\leq_{W_A\cap A}$ and $wW_A$ has the order $\leq_A$ induced from $W$.
\end{lemma}
\begin{proof}
    Let $t$ be a reflection of $W_A$. By Proposition \ref{cox:length:prop}(\ref{cox:length:prop:mult-char}), one has
    \begin{equation*}
        \ell_A'(tv^{-1})<\ell_A'(v^{-1})\iff t\in N_A(v^{-1})+v^{-1}(A\cap W_A)v
    \end{equation*}
    for all $v\in W_A$. But
    \begin{align*}
        N_A(v^{-1})&=A + v^{-1}A,\\
        &=A + N(v^{-1}) + v^{-1}A v,\\
        &=(A + N(v^{-1}) + v^{-1}A v)\cap W_A,\\
        &=(A\cap W_A) + (N(v^{-1})\cap W_A) + v^{-1}(A\cap W_A)v.
    \end{align*}
    Hence
    \begin{equation*}
        \ell_A'(tv^{-1})<\ell_A'(v^{-1})\iff t\in (A\cap W_A) + (N(v^{-1})\cap W_A).
    \end{equation*}
    By definition, $u^{-1}\in X_{W_A}$ implies $N(u^{-1})\cap W_A=\emptyset$ and thus $v^{-1}N(u^{-1})v\cap W_A=\emptyset$. It follows that 
    \begin{align*}
        N(v^{-1})\cap W_A&=(N(v^{-1}) + v^{-1}N(u^{-1})v)\cap W_A,\\
        &=N(v^{-1}u^{-1})\cap W_A.
    \end{align*}
    Altogether, using the fact that $t\in W_A$ and $w=uv$, one has
    \begin{align*}
        \ell_A'(tv^{-1})<\ell_A'(v^{-1})&\iff t\in A+ N(v^{-1}u^{-1}),\\
        &\iff \ell_A(wt)<\ell_A(w)
    \end{align*}
    where the final equivalence is an application of Proposition \ref{cox:length:prop}(\ref{cox:length:prop:mult-char}). The length additivity equality now follows from induction on $\ell'(v^{-1})$ using the same argument as Lemma \ref{cox:cosets:lem:poset-isom-right}. The second statement now follows.
\end{proof}

Everything is now in place to prove Proposition \ref{cox:cosets:prop:posets}.

\begin{proof}[Proof of Proposition \ref{cox:cosets:prop:posets}]
    For (\ref{cox:cosets:prop:posets:left}) note that $vA\cap W_J=\emptyset$ by definition of $^JW_A$. The result then follows from Lemma \ref{cox:cosets:lem:poset-isom-right}. The proof of (\ref{cox:cosets:prop:posets:right}) is more involved. Write $u=xy^{-1}$ where $y\in W_A$ and $x^{-1}\in X_{W_A}$. Then $N(u^{-1})\cap W_A=N(y)\cap W_A$. Observe that
    \begin{align*}
        N_A(y)&=A+yA,\\
        &=A+N(y)+yA y^{-1},\\
        &=(A+N(y)+yA y^{-1})\cap W_A,\\
        &=(A\cap W_A)+(N(y)\cap W_A)+(yA y^{-1}\cap W_A),\\
        &=(A\cap W_A)+(N(y)\cap W_A)+y(A\cap W_A)y^{-1}.
    \end{align*}
    It follows that 
    \begin{align*}
        (N(u^{-1})+A)\cap W_A&=(N(y)+A)\cap W_A,\\
        &=(N(y)\cap W_A)+(A\cap W_A),\\
        &=N_A(y)+y(A\cap W_A)y^{-1},
    \end{align*}
    where the last expression is simply $y(A\cap W_A)$ in $W_A$. The hypothesis that $u\in W_A^K$ implies $(N(u^{-1})+A)\cap (W_A)_K=\emptyset$, hence also $y(A\cap W_A)\cap (W_A)_K=\emptyset$. This is analogous to the condition used in the proof of (\ref{cox:cosets:prop:posets:left}). To obtain the result, first write $uv=xy^{-1}v=x(v^{-1}y)^{-1}$. Then, apply Lemma \ref{cox:cosets:lem:poset-isom-left} to $w=uv$, followed by (\ref{cox:cosets:prop:posets:left}) in $W_A$, and then Lemma \ref{cox:cosets:lem:poset-isom-left} once more. In particular, one has
    \begin{align*}
        \ell_A(w)&=\ell_A(x)+\ell_A'(v^{-1}y),\\
        &=\ell_A(x)+\ell'(v^{-1})+\ell_A'(y),\\
        &=\ell_A(xy^{-1})+\ell'(v^{-1}),\\
        &=\ell_A(u)+\ell'(v),
    \end{align*}
    as required. The final equality makes use of the fact that $\ell'(v^{-1})=\ell'(v)$ for the standard length function $\ell'$ of $W_A$.
\end{proof}

The next two results prove the generalisations of Kilmoyer and Howlett, respectively.

\begin{proof}[Proof of Proposition \ref{ctd:prelim:cosets:prop:kilmoyer}]
    The inclusion $W_L\subseteq W_J\cap w(W_A)_Kw^{-1}$ is immediate. For the converse, let $u\in W_J\cap w(W_A)_Kw^{-1}$. Thus $uw=wv$ for some $v\in (W_A)_K$. Applying both sides to $A$ implies that
    \begin{equation*}
        N(u)+u(wA)u^{-1}=wA+wN_A(v)w^{-1}.
    \end{equation*}
    As $u\in W_J$ and $wA\cap W_J$ is empty, intersecting both sides with $W_J$ gives
    \begin{equation*}
        N(u)=wN_A(v)w^{-1}\cap W_J.
    \end{equation*}
    If $u\neq1$ then 
    \begin{equation*}
        t=w^{-1}sw\in N_A(v)\subseteq (W_A)_K
    \end{equation*}
    for any left descent $s\in J$ of $u$. Applying both equalities of Proposition \ref{cox:cosets:prop:posets} yields
    \begin{equation*}
        \ell_A(w)+1=\ell_A(sw)=\ell_A(wt)=\ell_A(w)+\ell'(t).
    \end{equation*}
    In particular, $t\in K$ and so $s\in L$. But now $su\in W_J\cap w(W_A)_Kw^{-1}$ and $\ell(su)<\ell(u)$. By induction on $\ell$ one concludes that $su\in W_L$. Hence $u\in W_L$.
\end{proof}

\begin{proof}[Proof of Proposition \ref{ctd:prelim:cosets:prop:howlett}]
    Let $w'\in W_Jw(W_A)_K$, so $w'=a'wc'$ for some $a'\in W_J$ and $c'\in (W_A)_K$. As an element of $(W_A)_K$, one can write $c'=c''c$ where $c''\in (W_A)_{K\cap w^{-1}Jw}$ and $c\in (W_A)_K\cap{^{K\cap w^{-1}Jw}}(W_A)_\emptyset$. Then 
    \begin{equation*}
        a'wc'=a'(wc''w^{-1})wc
    \end{equation*}
    and moreover $a''=a'(wc''w^{-1})\in W_J$. If $a''wc=a''_1wc_1$ for some $a_1''\in W_J$ and $c_1\in (W_A)_K\cap{^{K\cap w^{-1}Jw}}(W_A)_\emptyset$ then $a''=a''_1$ and $c=c_1$. Indeed, since $c,c_1\in (W_A)_K$ and $a'',a''_1\in W_J$, one has
    \begin{equation*}
        cc_1^{-1}=w^{-1}((a'')^{-1}a''_1)w\in (W_A)_K\cap w^{-1}W_Jw=(W_A)_{K\cap w^{-1}Jw}
    \end{equation*}
    where the final equality is a result of Proposition \ref{ctd:prelim:cosets:prop:kilmoyer}. But $c$ and $c_1$ are unique minimal right coset representatives for $(W_A)_{K\cap w^{-1}Jw}$ in $(W_A)_K$. Then $cc_1^{-1}\in (W_A)_{K\cap w^{-1}Jw}$ implies that $c$ and $c_1$ are contained in the same right coset of $(W_A)_{K\cap w^{-1}Jw}$. It must therefore follow that $c=c_1$, then also $a''=a''_1$. To complete the uniqueness of expression statement, notice that the decomposition $a''=ab$ is unique, where $a\in W_J\cap W^{J\cap wKw^{-1}}$ and $b\in W_{J\cap wKw^{-1}}$. Now, observe that 
    \begin{align*}
        wcA\cap W_J&=(N(w)+w(cA)w^{-1})\cap W_J,\\
        &=(N(w)+w(N_A(c)+A)w^{-1})\cap W_J,\\
        &=(N(w)+wN_A(c)w^{-1}+wA w^{-1})\cap W_J,\\
        &=(wA+wN_A(c)w^{-1})\cap W_J,\\
        &=(wA\cap W_J)+(wN_A(c)w^{-1}\cap W_J).
    \end{align*}
    Since $w\in {^JW_A}$ one has $wA\cap W_J=\emptyset$. Moreover, $N_A(c)\subseteq (W_A)_K$ because $c\in (W_A)_K$. This implies that
    \begin{equation*}
        N_A(c)\cap w^{-1}W_Jw=N_A(c)\cap (W_A)_{K\cap w^{-1}Jw}=\emptyset
    \end{equation*}
    where the first equality is an application of Proposition \ref{ctd:prelim:cosets:prop:kilmoyer} and the second equality follows from the fact that $c\in {^{K\cap w^{-1}Jw}}(W_A)_\emptyset$. Altogether, this implies $wcA\cap W_J=\emptyset$ and thus $wc\in{^JW_A}$. Then, applying Proposition \ref{cox:cosets:prop:posets} three times, yields
    \begin{align*}
        \ell_A(a''wc)&=\ell(a'')+\ell_A(wc),\\
        &=\ell(a)+\ell(b)+\ell_A(wc),\\
        &=\ell(a)+\ell(b)+\ell_A(w)+\ell'(c).
    \end{align*}
    Observe that $abwc=aw(w^{-1}bw)c$ where $b'=w^{-1}bw\in (W_A)_{K\cap w^{-1}Jw}$. By Proposition \ref{ctd:prelim:cosets:prop:kilmoyer}, conjugation by $w^{-1}$ induces a length-preserving isomorphism $W_{wKw^{-1}\cap J}\simeq (W_A)_{K\cap w^{-1}Jw}$ of Coxeter groups, whence $\ell(b)=\ell'(b')$. Finally, notice that 
    \begin{align*}
        \ell_A(awb'c)&=\ell_A(abwc),\\
        &=\ell(a)+\ell(b)+\ell_A(w)+\ell'(c),\\
        &=\ell(a)+\ell'(b')+\ell_A(w)+\ell'(c),\\
        &=\ell_A(aw)+\ell'(b'c),
    \end{align*}
    which, because $b'c$ ranges over all elements of $(W_A)_K$, implies that $aw\in W_A^K$.
\end{proof}


\end{document}